\documentclass[a4paper,11pt]{article}
\usepackage[all,ps,arc]{xy}  
\usepackage{amsmath,amssymb,amsthm,latexsym}
\usepackage{pdfsync}
\usepackage{makeidx}
\usepackage[english]{babel}
\usepackage[applemac]{inputenc}
\usepackage{graphicx}
\usepackage[x11names]{xcolor}
\usepackage{geometry}
\newtheorem{Theorem}{Theorem}[section]
\newtheorem{Lemma}[Theorem]{Lemma}
\newtheorem{Proposition}[Theorem]{Proposition}
\newtheorem{Definition}[Theorem]{Definition}
\newtheorem{Corollary}[Theorem]{Corollary}
\newtheorem{Remark}[Theorem]{Remark}
\newtheorem{Example}[Theorem]{Example}

\newtheorem{Condition}[Theorem]{Condition}
\newtheorem{Text}[Theorem]{}

\newcommand{\cA}{{\mathcal A}}
\newcommand{\cB}{{\mathcal B}}
\newcommand{\cC}{{\mathcal C}}

\newcommand{\cF}{{\mathcal F}}

\newcommand{\cN}{{\mathcal N}}

\newcommand{\cS}{{\mathcal S}}

\newcommand\Arr{\mathbf{Arr}}
\newcommand\Seq{\mathbf{Seq}}
\newcommand\IsoSeq{\mathbf{IsoSeq}}
\newcommand\Ch{\mathbf{Ch}}

\newcommand\Grpd{\mathbf{Grpd}}

\newcommand\Ker{\mathrm{Ker}}
\newcommand\Cok{\mathrm{Cok}}
\newcommand\Dom{\mathrm{Dom}}
\newcommand\Cod{\mathrm{Cod}}

\newcommand\id{\mathrm{id}}

\begin{document}

\title{The snail lemma and the long homology sequence}

\author{Julia Ramos Gonz\'alez\footnote{Educational Centre for Mathematics, Education, Econometrics and Statistics, 
KULeuven, Warmoesberg 26, 1000 Brussel, Belgi$\ddot{\mbox{e}}$, julia.ramosgonzalez@kuleuven.be}
\footnote{Most of the research presented in this work was carried out while the first-named author was a Postdoctoral 
Researcher of the Fonds de la Recherche Scientifique - FNRS grant 32709538} 
\;and Enrico M. Vitale\footnote{Institut de recherche en math\'ematique et physique, UCLouvain,
Chemin du Cyclotron 2, 1348 Louvain-la-Neuve, Belgique, enrico.vitale@uclouvain.be}} 

\maketitle

\noindent {\bf Abstract:}
\noindent In the first part of the paper, we establish an homotopical version of the snail lemma (which is a generalization of 
the classical snake lemma). In the second part, we introduce the category $\Seq(\cA)$ of sequentiable families of arrows in 
a category $\cA$ and we compare it with the category of chain complexes in $\cA.$ We apply the homotopy snail lemma to 
a morphism in $\Seq(\cA)$ obtaining first a six-term exact sequence in $\Seq(\cA)$ and then, unrolling the sequence in 
$\Seq(\cA),$ a long exact sequence in $\cA.$ When $\cA$ is abelian, this sequence subsumes the usual long homology 
sequence obtained from an extension of chain complexes. 
\\ \\
\noindent{\it Keywords:} homology, exact sequence, homotopy kernel, snail lemma, sequentiable family. \\
{\it 2020 MSC:} 18G50, 18G35, 18N40.

\tableofcontents

\section{Introduction}\label{SecIntro}

A cornerstone in homological algebra is the fact that, starting from a short exact sequence of chain complexes in an 
abelian category $\cA,$ one can construct a long exact sequence relating the homology objects of the original complexes. 
A classical strategy to prove such a theorem is to prove first the snake lemma and then to construct the long exact sequence 
in homology pasting together the infinitely many six-term exact sequences coming from the snake lemma (see for 
example \cite{WBook}).

The snake lemma is a special case of a more general result, the snail lemma, introduced in \cite{EVSnail, ZJEVSnail} and exploited in
\cite{JMMVFract, MMVSnail} in order to unify some higher dimensional exact sequences appearing in homotopy theory, see \cite{GZ},
and in the study of groupoids and crossed modules, see \cite{RB, DKV}.The difference  between the snake and the snail lemmas lies in 
the fact that the snake requires as starting point a short exact sequence in $\Arr(\cA),$ the category of arrows in $\cA,$ whereas the 
snail works starting from any morphism in $\Arr(\cA).$ 

In this note, we show that it is possible to use the snail lemma instead of the snake lemma in order to construct a long exact 
sequence in homology starting from any morphism of chain complexes, and not necessarily from a short exact sequence of complexes. 
Even if the idea is quite simple, to state and prove it properly we have to introduce a new concept, that we call a sequentiable family of 
arrows. The idea behind a sequentiable family is to focus our attention not on the homology objects associated with a chain complex, 
but on the homology arrows, that is, those arrows whose kernel and cokernel are the homology objects associated with a complex
(see Section \ref{SecFromChToSeqFam}, and especially \ref{TextInfSeq}, for a more precise explication on how sequentiable families 
of arrows naturally arise from chain complexes). This allows us to formulate the snail lemma inside the category of sequentiable families, 
which is equipped with a structure of nullhomotopies more convenient than the one usually considered in the category of chain complexes.
This produces a single six-term exact sequence of sequentiable families, sequence which provides a compact presentation of a long exact 
sequence in the base category $\cA$ (Corollary \ref{CorLongSeq}) and, as a special case, the long homology sequence (Corollary 
\ref{CorLongExSeq}).

The layout of the paper is as follows: in Sections \ref{SecNullHom}, \ref{SecHomSnailSeq} and \ref{SecExactSnail}, we investigate a
general version of the snail lemma in a pointed category equipped with a structure of nullhomotopies, the guiding example being 
$\Arr(\cA)$ for $\cA$ an abelian category. We get a six-term sequence which is exact relatively to a good class of ``surjections''.
In the second part of the paper, Sections \ref{SecSeqFam} and \ref{SecFromChToSeqFam}, we introduce the new category $\Seq(\cA)$ 
of sequentiable families and we compare its nullhomotopy structure with that of the category $\Ch(\cA)$ of chain complexes when $\cA$
is preadditive. Finally, we apply the homotopy snail lemma in $\Seq(\cA)$ in order to get first a six-term exact sequence in $\Seq(\cA)$ 
and then, unrolling the sequence in $\Seq(\cA),$ a long exact sequence in $\cA.$ If the morphism in $\Seq(\cA)$ 
that we use to construct the long exact sequence comes from an exension of complexes, the new sequence coincides with the classical long 
homology sequence. Starting from Section \ref{SecExactSnail}, regular, protomodular, preadditive and abelian categories will appear so to 
make the base category $\cA$ reach enough. Basic references for these kinds of categories are \cite{BO2,BB,BG}.


\noindent N.B.: Given two arrows $\xymatrix{A \ar[r]^-{f} & B \ar[r]^-{g} & C},$ the composite arrow will be written as $f \cdot g.$

\section{Preliminaries on nullhomotopies}\label{SecNullHom} 

Categories with a structure of nullhomotopies have been introduced by Grandis in \cite{GR01}. 
We follow \cite{MMV23, VIT23} for the definition of nullhomotopy structure and of strong homotopy (co)kernel.

\begin{Definition}\label{DefNullHom}{\rm
A structure of nullhomotopies $\Theta$ on a category $\cB$ is given by:
\begin{enumerate}
\item[1)] For every arrow $g$ in $\cB,$ a set $\Theta(g)$ whose elements are called nullhomotopies on $g.$
\item[2)] For every triple of composable arrows $\xymatrix{A \ar[r]^f & B \ar[r]^g & C \ar[r]^h & D},$ a map 
$$f \circ - \circ h \colon \Theta(g) \to \Theta(f \cdot g \cdot h)$$
in such a way that, for every $\varphi \in \Theta(g),$ one has
\begin{enumerate}
\item $(f' \cdot f) \circ \varphi \circ (h \cdot h') = f' \circ (f \circ \varphi \circ h) \circ h'$ whenever the compositions $f' \cdot f$ and $h \cdot h'$ 
are defined,
\item $\id_B \circ \varphi \circ \id_C = \varphi.$
\end{enumerate}
\end{enumerate}
Equivalently, $\Theta$ is a functor from the twisted arrow category to the category of sets. See Remark 2.1.10 in \cite{VIT24} and the references 
therein for more details. \\
When $f=\id_B$ or $h=\id_C,$ we write $\varphi \circ h$ and $f \circ \varphi$ instead of $\id_B \circ \varphi \circ h$ 
and $f \circ \varphi \circ \id_C.$ \\
We sometimes depict a nullhomotopy $\lambda \in \Theta(g)$ as
$$\xymatrix{A \ar@/^1.0pc/[rr]^{g} \ar@/_1.0pc/@{--}[rr]_{} \ar@{}[rr]|{\lambda \; \Uparrow} & & B}$$
}\end{Definition}

Here is an additional condition on a structure of nullhomotopies introduced in \cite{GR01}.

\begin{Condition}\label{CondRedInter}{\rm
Let $(\cB,\Theta)$ be a category with nullhomotopies. 
The structure $\Theta$ satisfies the reduced interchange if, for any pair of composable arrows $f \colon A \to B$ and $g \colon B \to C$
 and for any pair of nullhomotopies $\alpha \in \Theta(f)$ and $\beta \in \Theta (g),$ one has $\alpha \circ g = f \circ \beta.$ 
 }\end{Condition}
 
 In the next definition, we introduce homotopy kernels and homotopy cokernels. They are the only kinds of homotopy limits or colimits
 we need in this paper.
 
\begin{Definition}\label{DefHKer}{\rm
Let $g \colon B \to C$ be an arrow in a category with nullhomotopies $(\cB,\Theta).$ 
\begin{enumerate}
\item An homotopy kernel of $g$ with respect to $\Theta$ (or $\Theta$-kernel) is a universal triple 
$$\xymatrix{\cN(g) \ar[rr]_{n_g} \ar@/^2pc/@{--}[rrrr]^{}_{\Downarrow \; \nu_g} & & B \ar[rr]_{g} & & C }$$
This means that, for any other triple $(A \in \cB,f \colon A \to B,\varphi \in \Theta(f \cdot g)),$
there exists a unique arrow $f' \colon A \to \cN(g)$ such that $f' \cdot n_g = f$ and $f' \circ \nu_g = \varphi$
\item A $\Theta$-kernel $(\cN(g),n_g,\nu_g)$ is strong if, for any triple $(A,f \colon A \to \cN(g),\varphi \in \Theta(f \cdot n_g))$ such that 
$\varphi \circ g = f \circ \nu_g,$ there exists a unique nullhomotopy $\varphi' \in \Theta(f)$ such that $\varphi' \circ n_g = \varphi.$
\item The notion of (strong) $\Theta$-cokernel is dual of the notion of (strong) $\Theta$-kernel. The notation is
$$\cC(g) \in \cB, c_g \colon C \to \cC(g), \gamma_g \in \Theta(g \cdot c_g)$$
\end{enumerate}
}\end{Definition}

\begin{Text}\label{TextCanc}{\rm
Let us recall from \cite{MMV23,VIT23} also the cancellation properties satisfied by a $\Theta$-kernel (the first one holds in general, 
the second one requires the reduced interchange):
\begin{enumerate}
\item Given arrows $g \colon B \to C$ and $f,h \colon A \to \cN(g),$ if $f \cdot n_g = h \cdot n_g$ and $f \circ \nu_g = h \circ \nu_g,$ then $f=h.$
\item Given arrows $g \colon B \to C$ and $f \colon A \to \cN(g)$ and nullhomotopies $\varphi,\psi \in \Theta(f)$ such that 
$\varphi \circ n_g = \psi \circ n_g,$ if the $\Theta$-kernel of $g$ is strong, then $\varphi = \psi.$
\end{enumerate}
}\end{Text} 

\begin{Text}\label{TextStrongZero}{\rm
In order to construct the snail sequence, we assume that the category $\cB$ has a zero object 0 which is $\Theta$-strong.
This means that, for every object $X \in \cB,$ the set of nullhomotopies $\Theta(0^X \colon X \to 0)$ of the terminal arrow is reduced to a single
element, denoted by $\ast^X,$ and the set of nullhomotopies $\Theta(0_X \colon 0 \to X)$ of the initial arrow is reduced to a single element, denoted 
by $\ast_X.$ \\
If we assume the reduced interchange, we get a canonical nullhomotopy 
$$\ast^X_Y \in \Theta(0^X_Y \colon X \to 0 \to Y)$$ 
given by $\ast^X \circ 0_Y$ or, equivalently, by $0^X \circ \ast_Y.$ Observe that in general $\Theta(0^X_Y)$ is not reduced to the element $\ast^X_Y,$
Nevertheless, for all arrows $f \colon W \to X, g \colon X \to Y, h \colon Y \to Z$ and for all nullhomotopies $\varphi \in \Theta(g)$ we have:
$$f \circ \ast^X_Y \circ h = \ast^W_Z \;\;\mbox{ and }\;\; 0^W_X \circ \varphi = \ast^W_Y \;\; \mbox{ and }\;\; \varphi \circ 0^Y_Z = \ast^X_Z$$
}\end{Text} 

\begin{Text}\label{TextN0Y}{\rm
The first relevant special case of $\Theta$-kernel is the one of the initial arrow $0_Y$ of an object $Y \in \cB.$ Its universal property can be restated as 
follows: for any object $X \in \cB$ and for any nullhomotopy $\varphi \in \Theta(0^X_Y),$ there exists a unique arrow $g \colon X \to \cN(0_Y)$ such that
$g \circ \nu_{0_Y} = \varphi$
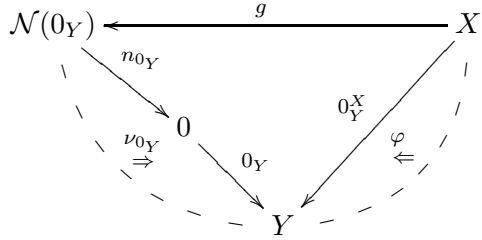
$$\xymatrix{\cN(0_Y) \ar[rd]^{n_{0_Y}} \ar@{}[rrdd]_{}_{\nu_{0_Y}}
& & & & X \ar[lldd]_{0^X_Y}^{\; \varphi} \ar[llll]_{g} \\
& 0 \ar[rd]^{0_Y} \\
& & Y \ar@/_2.0pc/@{--}[rruu]^{\; \Leftarrow} \ar@/^2pc/@{--}[lluu]_{\Rightarrow} }$$
}\end{Text} 

\begin{Text}\label{TextPi0}{\rm
Taking the $\Theta$-kernel of the arrow part of the $\Theta$-cokernel of the terminal arrow $0^Y,$ we get the $\pi_0$ of an object $Y.$
Here it is:
$$
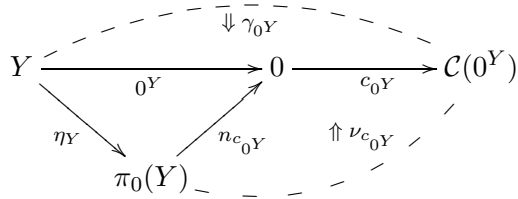
\xymatrix{Y \ar[rr]_{0^Y} \ar[rd]_{\eta_Y} \ar@/^2pc/@{--}[rrrr]^{}_{\Downarrow \; \gamma_{0^Y}} & & 
0 \ar[rr]_{c_{0^Y}} & & \cC(0^Y) \\
& \pi_0(Y) \ar[ru]_{n_{c_{0^Y}}} \ar@/_2pc/@{--}[rrru]^>>>>>>>>>>{\Uparrow \; \nu_{c_{0^Y}}} }$$
where $\eta_Y$ is the unique arrow such that $\eta_Y \circ \nu_{c_{0^Y}} = \gamma_{0^Y}.$
Observe that the construction of $\pi_0(Y)$ is a special case of the construction given in \ref{TextN0Y} since 
$\pi_0(Y) = \cN(c_{0^Y}) = \cN(0_{\cC(0^Y)}).$
}\end{Text} 

An interesting fact about objects of the form $\cN(0_Y)$ and, in particular, of the form $\pi_0(Y),$ is that they are discrete. Here is the 
general definition of discrete object. The name will be motivated at the end of Example \ref{ExArr1}.

\begin{Definition}\label{DefDiscr}{\rm
Consider a category with nullhomotopies $(\cB,\Theta).$ Assume that $\Theta$ satisfies the reduced interchange and that $\cB$ has 
a $\Theta$-strong zero object. An object $Y \in \cB$ is discrete if, for any arrow $g \colon X \to Y,$ the following conditions hold:
\begin{enumerate}
\item[-] if $g \neq 0^X_Y,$ then $\Theta(g) = \emptyset,$
\item[-] if $g = 0^X_Y,$ then $\Theta(g) = \{ \ast^X_Y \}.$
\end{enumerate}
Codiscrete objects are defined dually.
}\end{Definition} 


\begin{Lemma}\label{LemmaPi0Discr}
Consider a category with nullhomotopies $(\cB,\Theta).$ Assume that $\Theta$ satisfies the reduced interchange and that $\cB$ 
has a $\Theta$-strong zero object, strong $\Theta$-kernels and strong $\Theta$-cokernels.
For any object $Y \in \cB,$ the object $\cN(0_Y)$ is discrete and the object $\cC(0^Y)$ is codiscrete. In particular, $\pi_0(Y)$ is discrete.
\end{Lemma}

\begin{proof}
Consider an arrow $g \colon X \to \cN(0_Y)$ and a nullhomotopy $\varphi \in \Theta(g).$ We have to prove
that $g=0^X_{\cN(0_Y)}$ and $\varphi = \ast^X_{\cN(0_Y)}.$ Using the universal property of the $\Theta$-kernel, 
the first condition follows from the equations: 
\begin{enumerate}
\item[-] $g \cdot n_{0_Y} = 0^X = 0^X_{\cN(0_Y)} \cdot n_{0_Y}$
\item[-] $g \circ \nu_{0_Y} = \varphi \circ 0^{\cN(0_Y)}_Y = \ast^X_Y = \ast^X_{\cN(0_Y)} \circ 0^{\cN(0_Y)}_Y = 0^X_{\cN(0_Y)} \circ \nu_{0_Y}$
\end{enumerate} 
(in the second one we use \ref{TextStrongZero}). Using that the $\Theta$-kernel is strong, the second condition follows form the equation
\begin{enumerate}
\item[-] $\varphi \circ n_{0_Y} = \ast^X = \ast^X_{\cN(0_Y)} \circ n_{0_Y}$
\end{enumerate}
where we have used once again \ref{TextStrongZero}. The proof that $\cC(0^Y)$ is codiscrete is dual.
\end{proof} 

\begin{Lemma}\label{LemmaHomUsKer}
Consider a category with nullhomotopies $(\cB,\Theta).$ Assume that $\Theta$ satisfies the reduced interchange and that $\cB$ 
has a $\Theta$-strong zero object and $\Theta$-kernels.
If an object $Y \in \cB$ is discrete, then for every arrow $g \colon X \to Y$ the $\Theta$-kernel of $g$ coincide with
the usual categorical kernel of $g.$ 
\end{Lemma}

\begin{proof}
Since $Y$ is discrete, $n_g \cdot g = 0^{\cN(g)}_Y$ and $\nu_g = \ast^{\cN(g)}_Y.$ 
Consider an arrow $f \colon W \to X$ such that $f \cdot g = 0^W_Y.$ This implies that $\ast^W_Y \in \Theta(f \cdot g),$ 
so that there exists a unique arrow $f' \colon W \to \cN(g)$ such that $f' \cdot n_g = f$ and $f' \circ \nu_g = \ast^W_Y.$ 
If $f'' \colon W \to \cN(g)$ is another arrow such that $f'' \cdot n_g = f,$ then $f'' \circ \nu_g = f'' \circ \ast^{\cN(g)}_Y = \ast^W_Y$
by  \ref{TextStrongZero}, so that $f'' = f'.$ This proves that the $\Theta$-kernel satisfies the universal property of the usual kernel. \\
Conversely, let $k_g \colon \Ker(g) \to X$ be the usual kernel of $g.$ Since $k_g \cdot g = 0^{\Ker(g)}_Y,$ we can take
$\nu_g = \ast^{\Ker(g)}_Y \in \Theta(k_g \cdot g).$ Given now $f \colon W \to X$ and $\varphi \in \Theta(f \cdot g),$ since $Y$ is
discrete we get $f \cdot g = 0^W_Y$ and $\varphi = \ast^W_Y.$ From $f \cdot g = 0^W_Y,$ we get a unique $f' \colon W \to \Ker(g)$ 
such that $f' \cdot k_g = f.$ Moreover, $f' \circ \nu_g = f' \circ \ast^{\Ker(g)}_Y = \ast^W_Y = \varphi$ by \ref{TextStrongZero},
and we are done.
\end{proof} 

\begin{Text}\label{TextNId}{\rm
The second relevant special case of $\Theta$-kernel is the one of the identity arrow $\id_Y$ of an object 
$Y \in \cB.$ Its universal property can be restated as follows: for every arrow $g \colon X \to Y$ and for every 
nullhomotopy $\varphi \in \Theta(g),$ there exists a unique arrow $\overline{g} \colon X \to \cN(\id_Y)$ such that 
$\overline{g} \cdot n_{\id_Y} = g$ and $\overline{g} \circ \nu_{\id_Y} = \varphi$ 
$$\xymatrix{\cN(\id_Y) \ar[rrdd]^{n_{\id_Y}} \ar@{}[rrdd]_{}_{\nu_{\id_Y}}
& & & & X \ar[lldd]_{g}^{\; \varphi} \ar[llll]_{\overline{g}} \\ \\
& & Y \ar@/_2.0pc/@{--}[rruu]^{\; \Leftarrow} \ar@/^2pc/@{--}[lluu]_{\Rightarrow} }$$
Moreover, if the $\Theta$-kernel $\cN(\id_Y)$ is strong, there exists a unique nullhomotopy 
$\overline{\varphi} \in \Theta(\overline{g})$ such that $\overline{\varphi} \circ n_{\id_Y} = \varphi.$
}\end{Text} 

\begin{Text}\label{TextHomKerNId}{\rm
The $\Theta$-kernel of the identity arrow has a special role: if $\cB$ has pullbacks, then the $\Theta$-kernel of any 
arrow $g \colon X \to Y$ is given by
$$\xymatrix{\cN(g) \ar[rr]_-{n_g} \ar@/^2pc/@{--}[rrrr]^{}_{\Downarrow \; g' \circ \nu_{\id_Y}} & & X \ar[rr]_{g} & & Y }$$
where the diagram
$$\xymatrix{\cN(g) \ar[rr]^-{g'} \ar[d]_{n_g} & & \cN(\id_Y) \ar[d]^{n_{\id_Y}} \\
X \ar[rr]_{g} & & Y }$$
is a pullback (see Proposition 5.3 in \cite{MMV23}).
}\end{Text} 

The following general lemma is a variant of Remark 3.2 in \cite{VIT23}. Its dual also holds and both are needed in Section 
\ref{SecHomSnailSeq} to construct the snail sequence.

\begin{Lemma}\label{LemmaFunctHK}
Consider the solid part of the following commutative diagram (the one on the left) in a category with nullhomotopy $(\cB,\Theta)$
$$\xymatrix{\cN(a) \ar@{.>}[rr]^-{n(g,g')} \ar[d]^{n_a} \ar@/_3pc/@{--}[dd]^{ \; \nu_a \; \Rightarrow} 
& & \cN(b) \ar[d]_{n_b} \ar@/^3pc/@{--}[dd]_{\Leftarrow \; \nu_b \;}  \\
A \ar[rr]^-{g'} \ar[d]^{a} & & B \ar[d]_{b} \\
X \ar[rr]_-{g} & & Y }
\;\;\;\;\;\;\;\;\;\;\;\;
\xymatrix{\cN(a) \ar@{.>}[rr]_-{n(g,g')} \ar[d]^{n_a} \ar@/_2.5pc/@{--}[d]^{ \; \psi \; \Rightarrow} 
\ar@/^1.3pc/@{--}[rr]_-{n(\psi) \; \Downarrow} & & \cN(b) \ar[d]_{n_b}  \\
A \ar[rr]^-{g'} \ar[d]^{a} & & B \ar[d]_{b} \\
X \ar[rr]^-{g} \ar@/_1.3pc/@{--}[rr]^-{\varphi \; \Uparrow}  & & Y }$$
\begin{enumerate}
\item There exists a unique arrow $n(g,g') \colon \cN(a) \to \cN(b)$ such that $n(g,g') \cdot n_b = n_a \cdot g'$ and 
$n(g,g') \circ \nu_b = \nu_a \circ g.$
\item Consider the diagram on the right.
Given nullhomotopies $\varphi \in \Theta(g)$ and $\psi \in \Theta(n_a),$ if the $\Theta$-kernel of $b$ is strong
and if $\Theta$ satisfies the reduced interchange, then there exists a unique nullhomotopy $n(\psi) \in \Theta(n(g,g'))$ 
such that $n(\psi) \circ n_b = \psi \circ g'.$
\end{enumerate}
\end{Lemma}

\begin{proof}
1. Just apply the universal property of $\cN(b)$ to $\nu_a \circ g \in \Theta(n_a \cdot g' \cdot b).$ \\
2. Consider $\psi \circ g' \in \Theta(n(g,g') \cdot n_b).$ Since
$$\psi \circ g' \cdot b = \psi \circ a \cdot g = n_a \cdot a \circ \varphi = \nu_a \circ g = n(g,g') \circ \nu_b$$
we can use that $\cN(b)$ is strong to get a unique nullhomotopy $n(\psi) \in \Theta(n(g,g'))$ such that $n(\psi) \circ n_b = \psi \circ g'.$
\end{proof} 

\begin{Example}\label{ExArr1}{\rm
To help the reader with the various constructions introduced so far, let us look at the case where 
$(\cB,\Theta) = (\Arr(\cA),\Theta_{\Delta}).$ Objects, arrows and nullhomotopies for this example 
can be depicted as follows
$$\xymatrix{W \ar[rr]^-{f} \ar[d]_{w} & & X \ar[rr]^-{g} \ar[d]_{x} & & Y \ar[d]^{y} \ar[rr]^-{h} & & Z \ar[d]^{z} \\
W_0 \ar[rr]_-{f_0} & & X_0 \ar[rr]_-{g_0} \ar@{-->}[rru]^{\varphi} & & Y_0  \ar[rr]_-{h_0} & & Z_0 }$$
In other words, 
$$\Theta_{\Delta}((g,g_0) \colon (X,x,X_0) \to (Y,y,Y_0)) = \{ \varphi \colon X_0 \to Y \mid x \cdot \varphi = g, \varphi \cdot y = g_0 \}$$
and $(f,f_0) \circ \varphi \circ (h,h_0) = f_0 \cdot \varphi \cdot h.$
The structure $\Theta_{\Delta}$ satisfies the reduced interchange. Moreover, if $\cA$ has pullbacks, 
then $\Arr(\cA)$ has strong $\Theta_{\Delta}$-kernels constructed as in the following diagram on the left 
(the dashed arrow is the nullhomotopy), the one on the right being a pullback
$$\xymatrix{X \ar[rr]^-{\id_X} \ar[d]_{\langle x,g \rangle} & & X \ar[d]^>>>>{x} \ar[rr]^-{g} & & Y \ar[d]^{y} \\ 
X_0 \times_{g_0,y}Y \ar[rr]_-{y'} \ar@{-->}[rrrru]^<<<<<<<<{g_0'} & & X_0 \ar[rr]_-{g_0} & & Y_0 }
\;\;\;\;\;\; \xymatrix{X_0 \times_{g_0,y}Y \ar[rr]^-{g_0'} \ar[d]_{y'} & & Y \ar[d]^{y} \\ X_0 \ar[rr]_{g_0} & & Y_0 } $$
Dually, if $\cA$ has pushouts, then $\Arr(\cA)$ has
strong $\Theta_{\Delta}$-cokernels. Finally, if $\cA$ has a zero object 0, then the object $(0,\id_0,0)$
is a $\Theta_{\Delta}$-strong zero object in $\Arr(\cA).$ More details can be found in \cite{MMV23, VIT23}.\\
Assume now that $\cA$ has a zero object, kernels and cokernels. For an object $(Y,y,Y_0)$ in $\Arr(\cA),$
the $\Theta_{\Delta}$-kernel $\cN(0_{(Y,y,Y_0)})$ with its structural nullhomotopy is given by
$$\xymatrix{0 \ar[rr] \ar[d] & & Y \ar[d]^{y} \\
\Ker(y) \ar[rr]_{0} \ar@{-->}[rru]^{k_y} & & Y_0 }$$
Dually, the $\Theta_{\Delta}$-cokernel $\cC(0^{(Y,y,Y_0)})$ with its structural nullhomotopy is given by
$$\xymatrix{Y \ar[rr]^-{0} \ar[d]_{y} & & \Cok(y) \ar[d] \\
Y_0 \ar[rr] \ar@{-->}[rru]^{c_y} & & 0 }$$
so that the canonical arrow $\eta_{(Y,y,Y_0)} \colon (Y,y,Y_0) \to \pi_0(Y,y,Y_0)$ is given by
$$\xymatrix{Y \ar[rr] \ar[d]_{y} & & 0 \ar[d] \\
Y_0 \ar[rr]_-{c_y} & & \Cok(y) }$$
Finally, the $\Theta_{\Delta}$-kernel $\cN(\id_{(Y,y,Y_0)})$ with its structural nullhomotopy is given by
$$\xymatrix{Y \ar[d]_{\id_Y} \ar[rr]^-{\id_Y} & & Y \ar[d]^{y} \\
Y \ar[rr]_-{y} \ar@{-->}[rru]^{\id_Y} & & Y_0 }$$
Observe also that an object $(X,x,X_0) \in \Arr(\cA)$ is discrete if and only if $X=0$ (and then, necessarily, 
$x = 0_{X_0}).$ For the only if part, just consider the $\Theta_{\Delta}$-kernel of $\id_{(X,x,X_0)}.$
}\end{Example} 

\section{The homotopy snail sequence}\label{SecHomSnailSeq}

In this section, we work  in a  category with nullhomotopies $(\cB,\Theta)$ satisfying the reduced interchange 
condition \ref{CondRedInter}. We assume the existence of a $\Theta$-strong zero object as in \ref{TextStrongZero}. 
We also assume the existence of the strong $\Theta$-kernel of any arrow and of the strong $\Theta$-cokernel of 
any terminal arrow.  

\begin{Text}\label{TextConstrSnail}{\rm
Consider an arrow $g$ in $\cB$ together with its $\Theta$-kernel
$$\xymatrix{\cN(g) \ar[rr]_{n_g} \ar@/^2pc/@{--}[rrrr]_{\Downarrow \; \nu_g} & & X \ar[rr]_{g} & & Y }$$
We are going to construct the ``snail'' sequence connecting six discrete objects
$$\xymatrix{\cN(0_{\cN(g)}) \ar[r]^-{n(n_g)} & \cN(0_X) \ar[r]^-{n(g)} & \cN(0_Y) \ar[r]^-{\delta} & 
\pi_0(\cN(g)) \ar[r]^-{\pi_0(n_g)} & \pi_0(X) \ar[r]^-{\pi_0(g)} & \pi_0(Y) }$$
where each pair of consecutive arrows gives, by composition, the zero arrow. Its exactness will be discussed in 
Section \ref{SecExactSnail}. We split the construction of the sequence in five steps.
}\end{Text} 

\noindent {\bf Step 1:} By applying the first part of Lemma \ref{LemmaFunctHK} to the situations 
$$\xymatrix{\cN(0_{\cN(g)}) \ar@{.>}[rr]^-{n(n_g)} \ar[d]^{n_{0_{\cN(g)}}} \ar@/_4pc/@{--}[dd]^{\; \nu_{0_{\cN(g)}} \; \Rightarrow} 
& & \cN(0_X) \ar[d]_{n_{0_X}} \ar@/^4pc/@{--}[dd]_{\Leftarrow \; \nu_{0_X} \;}  \\
0 \ar[rr]^-{\id_0} \ar[d]^{0_{\cN(g)}} & & 0 \ar[d]_{0_X} \\
\cN(g) \ar[rr]_-{n_g} & & X }
\;\;\;\;\;\;\;\;\;
\xymatrix{\cN(0_X) \ar@{.>}[rr]^-{n(g)} \ar[d]^{n_{0_X}} \ar@/_4pc/@{--}[dd]^{\; \nu_{0_X} \; \Rightarrow} 
& & \cN(0_Y) \ar[d]_{n_{0_Y}} \ar@/^4pc/@{--}[dd]_{\Leftarrow \; \nu_{0_Y} \;}  \\
0 \ar[rr]^-{\id_0} \ar[d]^{0_X} & & 0 \ar[d]_{0_Y} \\
X \ar[rr]_-{g} & & Y }$$
we get the dotted arrows $n(n_g)$ and $n(g),$ which are unique with $n(n_g) \circ \nu_{0_X} = \nu_{0_{\cN(g)}} \circ n_g$ 
and $n(g) \circ \nu_{0_Y} = \nu_{0_X} \circ g.$ If we apply now the second part of Lemma \ref{LemmaFunctHK} to
$$\xymatrix{\cN(0_{\cN(g)}) \ar[rr]^-{n(n_g)} \ar[d]^{n_{0_{\cN(g)}}} \ar@/_5pc/@{--}[d]^{\; \ast^{\cN(0_{\cN(g)})} \; \Rightarrow} 
& & \cN(0_X) \ar[rr]^-{n(g)} & & \cN(0_Y) \ar[d]^{n_{0_Y}} \\
0 \ar[rrrr]^-{\id_0} \ar[d]^{0_{\cN(g)}} & & & & 0 \ar[d]^{0_Y} \\
\cN(g) \ar[rr]^-{n_g} \ar@{--}@/_1.5pc/[rrrr]^{\Uparrow \; \nu_g} & & X \ar[rr]^{g} & & Y \\ {} }$$
we get a nullhomotopy in $\Theta(n(n_g) \cdot n(g)).$ Since $\cN(0_Y)$ is discrete 
(see Lemma \ref{LemmaPi0Discr}), we can conclude that $n(n_g) \cdot n(g) = 0.$

\noindent {\bf Step 2:} Consider now the diagram
$$\xymatrix{X \ar[rrrr]^-{g} \ar[d]^{\eta_X} \ar@/_2pc/[dd]_{0^X} \ar@{--}@/_5pc/[ddd]^{\; \gamma_{0^X} \; \Rightarrow} 
& & & & Y \ar[d]_{\eta_Y} \ar@/^2pc/[dd]^{0^Y} \ar@{--}@/^5pc/[ddd]_{\Leftarrow \; \gamma_{0^Y} \;} \\
\pi_0(X) \ar@{.>}[rrrr]^-{\pi_0(g)} \ar[d]^{n_{c_{0^X}}} \ar@{--}@/^4pc/[dd] _{\Leftarrow \; \nu_{c_{0^X}} \;} 
& & & & \pi_0(Y) \ar[d]_{n_{c_{0^Y}}} \ar@{--}@/_4pc/[dd] ^{\; \nu_{c_{0^Y} \; \Rightarrow}} \\
0 \ar[d]^{c_{0^X}} & & & & 0 \ar[d]_{c_{0^Y}} \\
\cC(0^X) \ar@{.>}[rrrr]_-{c(g)} & & & & \cC(0^Y) }$$
Starting from $g \colon X \to Y$ and applying the dual of Lemma \ref{LemmaFunctHK}, we get a unique arrow $c(g)$ 
such that $\gamma_{0^X} \circ c(g) = g \circ \gamma_{0^Y}.$ If we start now from the arrow $c(g) \colon \cC(0^X) \to \cC(0^Y)$
and we work as in Step 1, we get a unique arrow $\pi_0(g)$ such that $\pi_0(g) \circ \nu_{c_{0^Y}} = \nu_{c_{0^X}} \circ c(g).$
Moreover, $g \cdot \eta_Y = \eta_X \cdot \pi_0(g).$ To see this, we compose with $\nu_{c_{0^Y}}$
$$\eta_X \cdot \pi_0(g) \circ \nu_{c_{0^Y}} = \eta_X \circ \nu_{c_{0^X}} \circ c(g) = \gamma_{0^X} \circ c(g) = g \circ \gamma_{0^Y} =
g \cdot \eta_Y \circ \nu_{c_{0^Y}}$$
and we can conclude using the first part of \ref{TextCanc}. \\
\noindent If we repeat the same argument starting from $n_g \colon \cN(g) \to X$ instead of $g \colon X \to Y,$ we get arrows 
$$c(n_g) \colon \cC(0^{\cN(g)}) \to \cC(0^X) \;\; \mbox{ and } \;\ \pi_0(n_g) \colon  \pi_0(\cN(g)) \to \pi_0(X)$$
unique with $\gamma_{0^{\cN(g)}} \circ c(n_g) = n_g \circ \gamma_{0^X}$ and $\pi_0(n_g) \circ \nu_{c_{0^X}} = \nu_{c_{0^{\cN(g)}}} \circ c(n_g).$
Moreover, by \ref{TextCanc} we get $n_g \cdot \eta_X = \eta_{\cN(g)} \cdot \pi_0(n_g)$ as above. \\
\noindent Finally, the nullhomotopy $\nu_g \in \Theta(n_g \cdot g)$ and the dual of the second 
part of Lemma \ref{LemmaFunctHK} allow us to prove that the composite arrow
$$\xymatrix{\cC(0^{\cN(g)}) \ar[rr]^-{c(n_g)} & & \cC(0^X) \ar[rr]^-{c(g)} & & \cC(0^Y) }$$
is the zero arrow (because $\cC(0^Y)$ is codiscrete). By the second part of Lemma \ref{LemmaFunctHK} again, we can conclude
that also the composite arrow
$$\xymatrix{\pi_0(\cN(g)) \ar[rr]^-{\pi_0(n_g)} & & \pi_0(X) \ar[rr]^-{\pi_0(g)} & & \pi_0(Y) }$$ 
is the zero arrow, as needed. 
 
 \noindent {\bf Step 3:} In order to construct the connecting arrow $\delta \colon \cN(0_Y) \to \pi_0(\cN(g)),$ consider the following diagram
 $$\xymatrix{ & & & & 
 \cN(0_Y) \ar[d]^{n_{0_Y}} \ar[lldd]_{0^{\cN(0_Y)}_X} \ar@{.>}@/_2pc/[lllldd]_{\Delta} \ar@{--}@/^3.5pc/[dd]_{\Leftarrow \; \nu_{0_Y} \; } \\ 
 & & & & 0 \ar[d]^{0_Y} \\ 
 \cN(g) \ar[rr]^-{n_g} \ar[d]_{\eta_{\cN(g)}} \ar@{--}@/_2pc/[rrrr]^{\nu_g \; \Uparrow} & & X \ar[rr]^-{g} & & Y \\ 
 \pi_0(\cN(g)) }$$
 Since $\nu_{0_Y} \in \Theta(n_{0_Y} \cdot 0_Y) = \Theta(0^{\cN(0_Y)}_X \cdot g),$ the universal property of the $\Theta$-kernel of $g$ gives 
 a unique arrow $\Delta \colon \cN(0_Y) \to \cN(g)$ such that $\Delta \cdot n_g = 0^{\cN(0_Y)}_X$ and $\Delta \circ \nu_g = \nu_{0_Y}.$ We put
 $$\delta = \Delta \cdot \eta_{\cN(g)} \colon \cN(0_Y) \to \cN(g) \to \pi_0(\cN(g))$$
 
 \noindent {\bf Step 4:} To check that the composite arrow
 $$\xymatrix{\cN(0_X) \ar[rr]^-{n(g)} & & \cN(0_Y) \ar[rr]^-{\delta} & & \pi_0(\cN(g)) }$$
 is the zero arrow, consider the diagram 
 $$\xymatrix{\cN(\id_{\cN(g)}) \ar[dd]^{n_{\id_{\cN(g)}}} \ar@{--}@/_3.5pc/[dd]^{\nu_{\id_{\cN(g)}} \Rightarrow} 
 & & & \cN(\id_X) \ar@{.>}[lll]_-{r_X} \ar@{.>}[llldd]_{s_X} \ar[rdd]^{n_{\id_X}} 
 \ar@{--}@/_2.5pc/[rdd]^{\nu_{\id_X} \Rightarrow} \ar@{--}@/^2pc/[llldd]_{\nu_{X,g} \Leftarrow}
 & & & \cN(0_X) \ar@{.>}[lll]_-{t_X} \ar[ld]_{n_{0_X}} \ar[rr]^-{n(g)} \ar@{--}@/^2.0pc/[lldd]_{\Leftarrow \; \nu_{0_X}}
 & & \cN(0_Y) \ar[d]^{n_{0_Y}} \ar@{--}@/_3pc/[dd]^{\; \nu_{0_Y} \Rightarrow} \\
 & & & & & 0 \ar[ld]_{0_X} & & & 0 \ar[d]^{0_Y} \\
 \cN(g) \ar[rrrr]^-{n_g} \ar@{--}@/_2pc/[rrrrrrrr]^-{\Uparrow \; \nu_g} \ar[d]_{\eta_{\cN(g)}} & & & & X \ar[rrrr]^-{g} & & & & Y \\
 \pi_0(\cN(g)) }$$
 By the universal property of $\cN(\id_X),$ we get a unique arrow $t_X$ such that $t_X \cdot n_{\id_X} = 0^{\cN(0_X)}_X$ and 
 $t_X \circ \nu_{\id_X} = \nu_{0_X}.$ By the universal property of $\cN(g),$ we get a unique arrow $s_X$ such that $s_X \cdot n_g = n_{\id_X}$ 
 and $s_X \circ \nu_g = \nu_{\id_X} \circ g.$ Since $\nu_{\id_X} \in \Theta(n_{\id_X}) = \Theta(s_X \cdot n_g)$ and 
 $\nu_{\id_X} \circ g = s_X \circ \nu_g,$ the fact that $\cN(g)$ is strong gives us a unique nullhomotopy $\nu_{X,g} \in \Theta(s_X)$ such that
 $\nu_{X,g} \circ n_g = \nu_{\id_X}.$ Finally, since $\nu_{X,g} \in \Theta(s_X),$ by the universal property of $\cN(\id_{\cN(g)})$ we get a unique
 arrow $r_X$ such that $r_X \cdot n_{\id_{\cN(g)}} = s_X$ and $r_X \circ \nu_{\id_{\cN(g)}} = \nu_{X,g}.$ Now we can check that
 $$n(g) \cdot \Delta = t_X \cdot r_X \cdot n_{\id_{\cN(g)}}$$
 using the first part of \ref{TextCanc}:
 \begin{enumerate}
 \item[-] $n(g) \cdot \Delta \cdot n_g = n(g) \cdot 0^{\cN(0_Y)}_X = 0^{\cN(0_X)}_X = t_X \cdot n_{\id_X} = 
 t_X \cdot s_X \cdot n_g = t_X \cdot r_X \cdot n_{\id_{\cN(g)}} \cdot n_g$
 \item[-] $n(g) \cdot \Delta \circ \nu_g = n(g) \circ \nu_{0_Y} = \nu_{0_X} \circ g = 
 t_X \circ \nu_{\id_X} \circ g = t_X \cdot s_X \circ \nu_g = t_X \cdot r_X \cdot n_{\id_{\cN(g)}} \circ \nu_g$
 \end{enumerate}
 Thanks to the previous equation, we have that 
 $$t_X \cdot r_X \circ \nu_{\id_{\cN(g)}} \circ \eta_{\cN(g)} \in \Theta(t_X \cdot r_X \cdot n_{\id_{\cN(g)}} \cdot \eta_{\cN(g)}) = 
 \Theta(n(g) \cdot \Delta \cdot \eta_{\cN(g)}) = \Theta(n(g) \cdot \delta)$$ 
 and, since $\pi_0(\cN(g))$ is discrete (see Lemma \ref{LemmaPi0Discr}), we can conclude that $n(g) \cdot \delta$ is the zero arrow.
 
 \noindent {\bf Step 5:} The fact that the composite arrow
 $$\xymatrix{\cN(0_Y) \ar[rr]^-{\delta} & & \pi_0(\cN(g)) \ar[rr]^-{\pi_0(n_g)} & & \pi_0(X) }$$
 is the zero arrow is obvious:
 $$\delta \cdot \pi_0(n_g) = \Delta \cdot \eta_{\cN(g)} \cdot \pi_0(n_g) = \Delta \cdot n_g \cdot \eta_X = 
 0^{\cN(0_Y)}_X \cdot \eta_X = 0^{\cN(0_Y)}_{\pi_0(X)}$$
 
 \begin{Example}\label{ExArr2}{\rm
 We go back to Example \ref{ExArr1}, where $(\cB,\Theta) = (\Arr(\cA),\Theta_{\Delta}).$ We assume that $\cA$
 has a zero object, pullbacks and cokernels. Starting from a $\Theta_{\Delta}$-kernel
 $$\xymatrix{X \ar[rr]^-{\id_X} \ar[d]_{\langle x,g \rangle} & & X \ar[d]^>>>>{x} \ar[rr]^-{g} & & Y \ar[d]^{y} \\ 
X_0 \times_{g_0,y}Y \ar[rr]_-{y'} \ar@{-->}[rrrru]^<<<<<<<<{g_0'} & & X_0 \ar[rr]_-{g_0} & & Y_0 }$$
 the associated snail sequence is
 $$\xymatrix{0 \ar[r] \ar[d] & 0 \ar[r] \ar[d] & 0 \ar[rr] \ar[d] & & 0 \ar[r] \ar[d] & 0 \ar[r] \ar[d] & 0 \ar[d] \\
 \Ker\langle x,g \rangle \ar[r] & \Ker(x) \ar[r] & \Ker(y) \ar[rr]^-{\delta_0} & & \Cok\langle x,g \rangle \ar[r] & \Cok(x) \ar[r] & \Cok(y) }$$ 
 where the unlabelled arrows are the obvious ones and $\delta_0$ is given by
 $$\delta_0 = \langle 0,k_y \rangle \cdot c_{\langle x,g \rangle} \colon 
 \Ker(y) \longrightarrow X_0 \times_{g_0,y}Y \longrightarrow \Cok\langle x,g \rangle$$
 so that the bottom line precisely is the snail sequence appearing in \cite{ZJEVSnail, EVSnail}.
 }\end{Example} 
 
 We end this section with a lemma about the arrow $\Delta \colon \cN(0_Y) \to \cN(g)$ defined in Step 3. This lemma, needed in the 
 next section, is a generalization of Lemma 3.2 in \cite{MMVSnail}.
 
 \begin{Lemma}\label{LemmaDelta} (With the previous notation.) The diagram
 $$\xymatrix{\cN(0_Y) \ar[rr]^-{\Delta} & & \cN(g) \ar[rr]^-{n_g} & & X }$$
 is a kernel (in the usual categorical sense).
 \end{Lemma}
 
 \begin{proof}
 The fact that $\Delta \cdot n_g = 0^{\cN(0_Y)}_X$ is part of the definition of $\Delta$ (see Step 3 above). 
 Consider now an arrow $a \colon A \to \cN(g)$ such that $a \cdot n_g = 0^A_X.$ This implies that
 $a \cdot n_g \cdot g = 0^A_Y$ ant then $a \circ \nu_g \in \Theta(0^A_Y).$ By the universal property of $\cN(0_Y),$ 
 we get a unique arrow $b \colon A \to \cN(0_Y)$ such that $b \circ \nu_{0_Y} = a \circ \nu_g.$ To check that $b \cdot \Delta = a$
 we use part 1 of \ref{TextCanc}:
 \begin{enumerate}
 \item[-] $b \cdot \Delta \cdot n_g = b \cdot 0^{\cN(0_Y)}_X = 0^A_X = a \cdot n_g$
 \item[-] $b \cdot \Delta \circ \nu_g = b \circ \nu_{0_Y} = a \circ \nu_g$
 \end{enumerate}
 As far as the uniqueness of $b$ is concerned, let $b' \colon A \to \cN(0_Y)$ be such that $b' \cdot \Delta = a.$ It follows that
 $b' \circ \nu_{0_Y} = b' \cdot \Delta \circ \nu_g = a \circ \nu_g,$ so that $b'=b$ by definition of $b.$
 \end{proof} 
 
 \section{Exactness of the snail sequence}\label{SecExactSnail} 
 
 We work under the same assumptions as in Section \ref{SecHomSnailSeq}: the nullhomotopy structure $\Theta$ satisfies 
 the reduced interchange condition \ref{CondRedInter}, the category $\cB$ is equipped with a $\Theta$-strong zero object 
 and has all the needed strong $\Theta$-kernels and strong $\Theta$-cokernels. We moreover assume that $\cB$ has pullbacks.
 
 \begin{Condition}\label{CondSurjLike}{\rm
 We fix a class of arrows $\cS$ in $\cB$ satisfying the following conditions:
 \begin{enumerate}
 \item $\cS$ is closed under composition,
 \item $\cS$ is stable under pullbacks,
 \item $\cS$ contains all the identities,
 \item $\cS$ has the left cancellation property: if a composite $f \cdot g$ is in $\cS,$ then $g$ is in $\cS.$
 \end{enumerate}
 Note that such a class $\cS$ contains all the isomorphisms. Note also that, if we ask that all monomorphisms in $\cS$
 are isomorphisms, we get the notion of surjection-like class of arrows discussed in \cite{PAJ24} (where condition 4 is 
 called strong right cancellation property).
 }\end{Condition}
 
 \begin{Definition}\label{DefSExact}{\rm
 Consider the diagram in $(\cB, \Theta)$
 $$\xymatrix{W \ar[rr]_{f} \ar@/^2pc/@{--}[rrrr]_{\Downarrow \; \varphi} & & X \ar[rr]_{g} & & Y }$$
 We say that $(f,\varphi,g)$ is $\cS$-exact if the unique factorization of $(f,\varphi)$ through the $\Theta$-kernel 
 $(n_g,\nu_g)$ of $g$ is in $\cS.$ \\
 Note that if $Y$ is discrete, to be $\cS$-exact means that the unique factorization of $f$ through the categorical
 kernel of $g$ is in $\cS.$
 }\end{Definition} 
 
 \begin{Definition}\label{DefSPropSGlob}{\rm
 Let $Y$ be an object in $(\cB,\Theta).$
 \begin{enumerate}
 \item $Y$ is $\cS$-proper if $\overline{y} \colon \cN(\id_Y) \to \cN(\eta_Y)$ is in $\cS,$ where $\overline{y}$ is the unique 
 arrow such that $\overline{y} \cdot n_{\eta_Y} = n_{\id_Y}$ and $\overline{y} \circ \nu_{\eta_Y} = \nu_{\id_Y} \circ \eta_Y$
 $$\xymatrix{\cN(\id_Y) \ar[rr]_{n_{\id_Y}} \ar[rd]_{\overline{y}} \ar@{--}@/^2pc/[rr]^{}_-{\Downarrow \; \nu_{\id_Y}} 
 & & Y \ar[rr]^{\eta_Y} & & \pi_0(Y) \\
 & \cN(\eta_Y) \ar[ru]^{n_{\eta_Y}} \ar@{--}@/_1.5pc/[rrru]^{\Uparrow \; \nu_{\eta_Y}} }$$
 \item$Y$ is $\cS$-global if $\eta_Y \colon Y \to \pi_0(Y)$ is in $\cS.$
 \end{enumerate}
 }\end{Definition} 

\begin{Example}\label{ExArr3}{\rm
When $(\cB,\Theta) = (\Arr(\cA),\Theta_{\Delta}),$ the factorization $\overline{y}$ of Definition \ref{DefSPropSGlob} is 
essentially the factorization of the arrow part $y$ of an object $(Y,y,Y_0)$ through the kernel of its cokernel
$$\xymatrix{Y \ar[d]_{\id_Y} \ar[rr]^-{\id_Y} \ar[rdd]^<<<<<<<<{\id_Y} & & Y \ar[d]^{y} \ar[rr] & & 0 \ar[d] \\
Y \ar[rr]^-{y} \ar[rdd]_{\overline{y}} & & Y_0 \ar[rr]_-{c_y} & & \Cok(y) \\
& Y \ar[d]_>>>>>>>{\overline{y}} \ar[ruu]^>>>>>>>>{\id_Y} \\
& \Ker(c_y) \ar[ruu]_{k_{c_y}} }$$
}\end{Example} 

\begin{Proposition}\label{PropExactSnail1}
Consider an arrow $g$ in $(\cB,\Theta)$ together with its $\Theta$-kernel
$$\xymatrix{\cN(g) \ar[rr]_{n_g} \ar@/^2pc/@{--}[rrrr]_{\Downarrow \; \nu_g} & & X \ar[rr]_{g} & & Y }$$
If $Y$ and $\cN(g)$ are $\cS$-proper and $X$ is $\cS$-global where $\cS$ is a class of morphisms in $\cB$ satisfying 
Condition \ref{CondSurjLike}, then the associated snail sequence
$$\xymatrix{\cN(0_{\cN(g)}) \ar[r]^-{n(n_g)} & \cN(0_X) \ar[r]^-{n(g)} & \cN(0_Y) \ar[r]^-{\delta} & 
\pi_0(\cN(g)) \ar[r]^-{\pi_0(n_g)} & \pi_0(X) \ar[r]^-{\pi_0(g)} & \pi_0(Y) }$$
is $\cS$-exact in $\cN(0_X), \cN(0_Y)$ and $\pi_0(X).$
\end{Proposition}

\begin{proof} We use the constructions and the notations introduced in Section \ref{SecHomSnailSeq}. 
In particular, for the arrows $t_X, r_X$ and $s_X$ and for the nullhomotopy $\nu_{X,g},$ see Step 4 of Section \ref{SecHomSnailSeq}. \\
{\bf Exactness in $\cN(0_X)$:} we are going to prove that the diagram
$$\xymatrix{\cN(0_{\cN(g)}) \ar[rr]^-{n(n_g)} & & \cN(0_X) \ar[rr]^-{n(g)} & & \cN(0_Y) }$$
constructed in Step 1 of Section \ref{SecHomSnailSeq}, is a kernel. This implies the $\cS$-exactness because $\cS$ 
contains the isomorphisms. 
Consider an arrow $a \colon A \to \cN(0_X)$ such that $a \cdot n(g) = 0^A_{\cN(0_Y)}.$ As a preliminary step, we check that 
$a \cdot t_X \cdot r_X \cdot n_{\id_{\cN(g)}} = 0^A_{\cN(g)}$ using the cancellation property \ref{TextCanc}:
\begin{enumerate}
\item[-] $t_X \cdot r_X \cdot n_{\id_{\cN(g)}} \cdot n_g = t_X \cdot s_X \cdot n_g = t_X \cdot n_{\id_X} = 0^{\cN(0_X)}_X = 
0^{\cN(0_X)}_{\cN(g)} \cdot n_g$ and then $a \cdot t_X \cdot r_X \cdot n_{\id_{\cN(g)}} \cdot n_g = 
a \cdot 0^{\cN(0_X)}_{\cN(g)} \cdot n_g = 0^A_{\cN(g)} \cdot n_g$
\item[-] $a \cdot t_X \cdot r_X \cdot n_{\id_{\cN(g)}} \circ \nu_g = a \cdot t_X \cdot s_X \circ \nu_g = 
a \cdot t_X \circ \nu_{\id_X} \circ g = a \circ \nu_{0_X} \circ g = a \cdot n(g) \circ \nu_{0_Y} = 
\ast^A_{\cN(0_Y)} \circ 0^{\cN(0_Y)}_Y = \ast^A_{\cN(0_Y)} \circ 0^{\cN(0_Y)}_{\cN(g)} \cdot n_g \cdot g = 0^A_{\cN(g)} \circ \nu_g$
\end{enumerate} 
By the universal property of $\cN(0_{\cN(g)}),$ there exists a unique arrow $a' \colon A \to \cN(0_{\cN(g)})$ such that
$a' \circ \nu_{0_{\cN(g)}} = a \cdot t_X \cdot r_X \circ \nu_{\id_{\cN(g)}}.$ Now we show that $a' \cdot n(n_g) = a$ using once again
\ref{TextCanc}:
\begin{enumerate}
\item[-] $a' \cdot n(n_g) \cdot n_{0_X} = 0^A = a \cdot n_{0_X}$
\item[-] $a' \cdot n(n_g) \circ \nu_{0_X} = a' \circ \nu_{0_{\cN(g)}} \circ n_g = a \cdot t_X \cdot r_X \circ \nu_{\id_{\cN(g)}} \circ n_g =
a \cdot t_X \circ \nu_{X,g} \circ n_g = a \cdot t_X \circ \nu_{\id_X} = a \circ \nu_{0_X}$
\end{enumerate}
As far as the uniqueness of the factorization is concerned, let $a'' \colon A \to \cN(0_{\cN(g)})$ be such that $a'' \cdot n(n_g) = a.$
Since clearly $a' \cdot n_{0_{\cN(g)}} = 0^A = a'' \cdot n_{0_{\cN(g)}},$ to have $a'' = a'$ it remains to show that
$a' \circ \nu_{0_{\cN(g)}} = a'' \circ \nu_{0_{\cN(g)}}.$ Since the $\Theta$-kernel $\cN(g)$ is strong, it is enough to compose with $n_g \colon$
\begin{enumerate}
\item[-] $a' \circ \nu_{0_{\cN(g)}} \circ n_g = a' \cdot n(n_g) \circ \nu_{0_X} =  a'' \cdot n(n_g) \circ \nu_{0_X} = a'' \circ \nu_{0_{\cN(g)}} \circ n_g$
\end{enumerate}

\noindent {\bf Exactness in $\cN(0_Y)$:} We are going to prove that the unique arrow $\sigma$ making commutative the following diagram 
is in $\cS \colon$
$$\xymatrix{\cN(\delta) \ar[rr]^-{n_{\delta}} & & \cN(0_Y) \ar[rr]^-{\delta} & & \pi_0(\cN(g)) \\
& & \cN(0_X) \ar[llu]^{\sigma} \ar[u]_{n(g)} }$$
For this, consider the factorization $\Delta'$ of $n_{\delta} \cdot \Delta$ through the (homotopy) kernel of $\eta_{\cN(g)} \colon$
$$\xymatrix{\cN(\delta) \ar[rr]^-{n_{\delta}} \ar[rrrrd]_{\Delta'} & & \cN(0_Y) \ar[rr]^-{\Delta} & & \cN(g) \ar[rr]^-{\eta_{\cN(g)}} & & \pi_0(\cN(g)) \\
& & & & \cN(\eta_{\cN(g)}) \ar[u]_{n_{\eta_{\cN(g)}}} }$$
Consider also the factorization $\overline{z}$ obtained when, in Definition \ref{DefSPropSGlob}.1, we start from the object $\cN(g) \colon$
$$\xymatrix{\cN(\id_{\cN(g)}) \ar[rr]^-{n_{\id_{\cN(g)}}} \ar[rd]_{\overline{z}} & & \cN(g) \\
& \cN(\eta_{\cN(g)}) \ar[ru]_{n_{\eta_{\cN(g)}}} }$$
Using the arrows $\sigma, \Delta'$ and $\overline{z},$ we can build up the diagram
$$\xymatrix{ \cN(\delta) \ar[rrrr]^-{\Delta'} & & & & \cN(\eta_{\cN(g)}) \\
 \cN(0_X) \ar[u]^{\sigma} \ar[rr]_-{t_X} & & \cN(\id_X) \ar[rr]_-{r_X} & & \cN(\id_{\cN(g)}) \ar[u]_{\overline{z}} }$$
To check its commutativity, we compose with $n_{\eta_{\cN(g)}},$ which is a monomorphism because $\pi_0(\cN(g))$ is discrete:
\begin{enumerate}
\item[-] $\sigma \cdot \Delta' \cdot n_{\eta_{\cN(g)}} = \sigma \cdot n_{\delta} \cdot \Delta = n(g) \cdot \Delta =
t_X \cdot r_X \cdot n_{\id_{\cN(g)}} = t_X \cdot r_X \cdot \overline{z} \cdot n_{\eta_{\cN(g)}}$
\end{enumerate} 
By assumption, $\cN(g)$ is $\cS$-proper, that is, $\overline{z} \in \cS.$ Since $\cS$ is stable under pullbacks, to prove that $\sigma \in \cS$
it remains to show that the previous square is a pullback. For this, consider two arrows
$$a \colon A \to \cN(\id_{\cN(g)}) \;\;\mbox{ and }\;\; b \colon A \to \cN(\delta)$$
such that $a \cdot \overline{z} = b \cdot \Delta'.$ This equality implies that
\begin{enumerate}
\item[] $a \cdot n_{\id_{\cN(g)}} \cdot n_g = a \cdot \overline{z} \cdot n_{\eta_{\cN(g)}} \cdot n_g = 
b \cdot \Delta' \cdot n_{\eta_{\cN(g)}} \cdot n_g =
b \cdot n_{\delta} \cdot \Delta \cdot n_g = b \cdot n_{\delta} \cdot 0^{\cN(0_Y)}_X = 0^A_X$
\end{enumerate}
so that we can consider the nullhomotopy $a \circ \nu_{\id_{\cN(g)}} \circ n_g \in \Theta(a \cdot n_{\id_{\cN(g)}} \cdot n_g) = \Theta(0^A_X).$
By the universal property of $\cN(0_X),$ we get a unique arrow $c \colon A \to \cN(0_X)$ such that 
$c \circ \nu_{0_X} = a \circ \nu_{\id_{\cN(g)}} \circ n_g.$ We have to verify that $c$ is a factorization of $a$ and $b.$ To check that 
$c \cdot \sigma = b,$ we compose with $n_{\delta},$ which is a monomorphism because $\pi_0(\cN(g))$ is discrete, and then we use the 
universal property of $\cN(0_Y) \colon$
\begin{enumerate}
\item[-] $c \cdot \sigma \cdot n_{\delta} \circ \nu_{0_Y} = c \cdot n(g) \circ \nu_{0_Y} = c \circ \nu_{0_X} \circ g = 
a \circ \nu_{\id_{\cN(g)}} \circ n_g \cdot g = a \cdot n_{\id_{\cN(g)}} \circ \nu_g = b \cdot n_{\delta} \cdot \Delta \circ \nu_g = 
b \cdot n_{\delta} \circ \nu_{0_Y}$
\end{enumerate} 
(the fourth equality follows from the rduced interchenge). To check that $c \cdot t_X \cdot r_X = a,$ we use \ref{TextCanc} and we compose 
a first time with $n_{\id_{\cN(g)}}$ and a second time with 
$\nu_{\id_{\cN(g)}}.$ In the first case we use \ref{TextCanc} again, and in the second case we use that the $\Theta$-kernel $\cN(g)$ is strong:
\begin{enumerate}
\item[-] $c \cdot t_X \cdot r_X \cdot n_{\id_{\cN(g)}} \cdot n_g = c \cdot t_X \cdot s_X \cdot n_g = c \cdot t_X \cdot n_{\id_X} = 
c \cdot 0^{\cN(0_X)}_X = 0^A_X = a \cdot n_{\id_{\cN(g)}} \cdot n_g$
\item[-] $c \cdot t_X \cdot r_X \cdot n_{\id_{\cN(g)}} \circ \nu_g = c \cdot t_X \cdot s_X \circ \nu_g = c \cdot t_X \circ \nu_{\id_X} \circ g =
c \circ \nu_{0_X} \circ g = a \circ \nu_{\id_{\cN(g)}} \circ n_g \cdot g = a \cdot n_{\id_{\cN(g)}} \circ \nu_g$
\item[-] $c \cdot t_X \cdot r_X \circ \nu_{\id_{\cN(g)}} \circ n_g = c \cdot t_X \circ \nu_{X,g} \circ n_g = c \cdot t_X \circ \nu_{\id_X} =
c \circ \nu_{0_X} = a \circ \nu_{\id_{\cN(g)}} \circ n_g$
\end{enumerate} 
As far as the uniqueness of the factorization $c$ is concerned, let $c' \colon A \to \cN(0_X)$ be such that $c' \cdot \sigma = b$ and 
$c' \cdot t_X \cdot r_X = a.$ To verify that $c' = c,$ we go back to the definition of $c \colon$
\begin{enumerate}
\item[-] $c' \circ \nu_{0_X} = c' \cdot t_X \circ \nu_{\id_X} = c' \cdot t_X \circ \nu_{X,g} \circ n_g = 
c' \cdot t_X \cdot r_X \circ \nu_{\id_{\cN(g)}} \circ n_g = a \circ \nu_{\id_{\cN(g)}} \circ n_g$
\end{enumerate}

\noindent {\bf Exactness in $\pi_0(X)$:} We are going to prove that the unique arrow $\sigma$ making commutative the following diagram 
is in $\cS \colon$
$$\xymatrix{\cN(\pi_0(g)) \ar[rr]^-{n_{\pi_0(g)}} & & \pi_0(X) \ar[rr]^-{\pi_0(g)} & & \pi_0(Y) \\
& & \cN(0_X) \ar[llu]^{\sigma} \ar[u]_{\pi_0(n_g)} }$$
We can split the pullback describing the $\Theta$-kernel $\cN(g)$ (see \ref{TextHomKerNId}) in two steps:
$$\xymatrix{\cN(g) \ar[rrr]^-{g'} \ar[rd]^{a} \ar[dd]_{n_g} & & & \cN(\id_Y) \ar[dd]_<<<<<<{n_{\id_Y}} \ar[rd]^{\overline{y}} \\
& P \ar[ld]^{b} \ar[rrr]^<<<<<<<<<<<<{c} & & & \cN(\eta_Y) \ar[ld]^{n_{\eta_Y}} \\
X \ar[rrr]_-{g} & & &  Y }$$
Now observe that
\begin{enumerate}
\item[] $b \cdot \eta_X \cdot \pi_0(g) = b \cdot g \cdot \eta_Y = c \cdot n_{\eta_Y} \cdot \eta_Y = 
c \cdot 0^{\cN(\eta_Y)}_{\pi_0(Y)} = 0^{P}_{\pi_0(Y)}$
\end{enumerate} 
so that there exists a unique arrow $t \colon P \to \cN(\pi_0(g))$ such that $t \cdot n_{\pi_0(g)} = b \cdot \eta_X.$
In fact, more is true: the square
$$\xymatrix{P \ar[rr]^{t} \ar[d]_{b} & & \cN(\pi_0(g)) \ar[d]^{n_{\pi_0(g)}} \\
X \ar[rr]_-{\eta_X} & & \pi_0(X) }$$
is a pullback. This can be proved using the following commutative diagrams:
$$\xymatrix{P \ar[rr]^-{c} \ar[d]_{b} & & \cN(\eta_Y) \ar[d]^{n_{\eta_Y}} \ar[rr] & & 0 \ar[d] \\
X \ar[rr]_-{g} \ar@{}[rru]|{(1)} & & Y \ar[rr]_-{\eta_Y}  \ar@{}[rru]|{(2)} & & \pi_0(Y) }$$
$$\xymatrix{P \ar[rr]^-{t} \ar[d]_{b} & & \cN(\pi_0(g)) \ar[d]^{n_{\pi_0(g)}} \ar[rr] & & 0 \ar[d] \\
X \ar[rr]_-{\eta_X}  \ar@{}[rru]|{(3)} & & \pi_0(X) \ar[rr]_-{\pi_0(g)}  \ar@{}[rru]|{(4)} & & \pi_0(Y) }$$
Since (1) and (2) are pullbacks, then (1)+(2) is a pullback. This implies that (3)+(4) is a pullback because $\eta_X \cdot \pi_0(g) = g \cdot \eta_Y.$
Since (4) also is a pullback, we can conclude that (3) is a pullback. Finally, observe that $a \cdot t = \eta_{\cN(g)} \cdot \sigma.$ To check this
equality, we can compose with $n_{\pi_0(g)},$ which is a monomorphism because $\pi_0(Y)$ is discrete:
\begin{enumerate}
\item[-] $a \cdot t \cdot n_{\pi_0(g)} = a \cdot b \cdot \eta_X = n_g \cdot \eta_X = \eta_{\cN(g)} \cdot \pi_0(n_g) = 
\eta_{\cN(g)} \cdot \sigma \cdot n_{\pi_0(g)}$
\end{enumerate}
Now we can conclude as follows: since $Y$ is $\cS$-proper, $\overline{y} \in \cS$ and then $a \in \cS$ by stability under pullbacks. Since $X$
is $\cS$-global, $\eta_X \in \cS$ and then $t \in \cS$ once again by stability under pullbacks. Since $\cS$ is closed under composition, the equality
$a \cdot t = \eta_{\cN(g)} \cdot \sigma$ and the left cancellation property imply that $\sigma \in \cS.$
\end{proof}

The exactness of the snail sequence in $\pi_0(\cN(g))$ requires one more assumption on the class $\cS.$

\begin{Condition}\label{CondB}{\rm
Let $\cS$ be a class of arrows in $\cB$ as in Condition \ref{CondSurjLike}. We say that $\cS$ satisfies condition 
(Sub) if, in the following commutative diagram, if $K(x,y)$ and $g$ are in $\cS,$ then $K(g,g_0)$ also is in $\cS.$
$$\xymatrix{ & & \Ker(x) \ar[rr]^-{K(g,g_0)} \ar[d]_{k_x} & & \Ker(y) \ar[d]^{k_y} \\
\Ker(g) \ar[rr]^-{k_g} \ar[d]_{K(x,y)} & & X \ar[rr]^-{g} \ar[d]_{x} & & Y \ar[d]^{y} \\
\Ker(g_0) \ar[rr]_-{k_{g_0}} & & X_0 \ar[rr]_{g_0} & & Y_0 }$$
}\end{Condition} 

\begin{Remark}\label{RemZurab}{\rm
Condition (Sub), with $\cS = \{\mbox{regular epimorphisms}\},$ has been isolated by Bourn in \cite{DB3x3} (see also 
\cite{BB}) and is a special case of the snake lemma. It has been used in \cite{EVSnail} to prove the basic version of 
the snail lemma which, in the context of pointed protomodular regular categories, subsumes the snake lemma. The 
precise situation has been explained to us by Zurab Janelidze in a private communication and we report it here for 
the sake of completeness. Assume that the category $\cB$ is pointed and regular. Then condition (Sub) is equivalent 
to the subtractivity of $\cB.$ First, by Theorem 3 in \cite{ZJEVSnail}, the subtractivity of $\cB$ is equivalent to the fact 
that the incomplete snail lemma holds in $\cB.$ Now, the incomplete snail lemma is equivalent to condition (Sub). 
Indeed, if in the first diagram of Section 3 in \cite{ZJEVSnail} we assume that the arrows $Y_1 \to X$ and $W_1 \to X$ 
are monos, then the incomplete snail lemma precisely gives condition (Sub). For the converse implication, one has to 
consider the (regular epi, mono) factorization of the same two arrows. The fact that subtractivity implies condition (Sub) 
can also be checked using the pointed subobject functor introduced in \cite{ZJ}.
}\end{Remark} 

\begin{Proposition}\label{PropExactSnail2}
Consider an arrow $g$ in $(\cB,\Theta)$ together with its $\Theta$-kernel
$$\xymatrix{\cN(g) \ar[rr]_{n_g} \ar@/^2pc/@{--}[rrrr]_{\Downarrow \; \nu_g} & & X \ar[rr]_{g} & & Y }$$
If $X$ is $\cS$-proper and $\cN(g)$ is $\cS$-global and if $\cS$ satisfies Conditions \ref{CondSurjLike} and \ref{CondB}, 
then the associated snail sequence
$$\xymatrix{\cN(0_{\cN(g)}) \ar[r]^-{n(n_g)} & \cN(0_X) \ar[r]^-{n(g)} & \cN(0_Y) \ar[r]^-{\delta} & 
\pi_0(\cN(g)) \ar[r]^-{\pi_0(n_g)} & \pi_0(X) \ar[r]^-{\pi_0(g)} & \pi_0(Y) }$$
is $\cS$-exact in $\pi_0(\cN(g)).$
\end{Proposition}

\begin{proof}

We are going to prove that the unique arrow $\sigma$ making commutative the following diagram 
is in $\cS \colon$
$$\xymatrix{\cN(\pi_0(n_g)) \ar[rr]^-{n_{\pi_0(n_g)}} & & \pi_0(\cN(g)) \ar[rr]^-{\pi_0(n_g)} & & \pi_0(X) \\
& & \cN(0_Y) \ar[llu]^{\sigma} \ar[u]_{\delta} }$$
Consider the following diagram, where the square is the pullback describing the $\Theta$-kernel of $g$ (see \ref{TextHomKerNId})
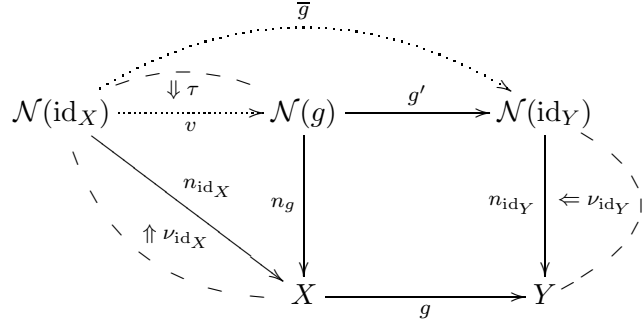
$$\xymatrix{\cN(\id_X) \ar@{.>}[rr]_-{v} \ar[rrdd]^{n_{\id_X}} \ar@{.>}@/^2.7pc/[rrrr]^-{\overline{g}}
\ar@{--}@/_2pc/[rrdd]^{\Uparrow \; \nu_{\id_X}} \ar@{--}@/^1.3pc/[rr]_{\Downarrow \; \tau}
& & \cN(g) \ar[rr]^-{g'} \ar[dd]_{n_g} & & \cN(\id_Y) \ar[dd]_{n_{\id_Y}} \ar@{--}@/^3pc/[dd]_{\Leftarrow \; \nu_{\id_Y}} \\ \\
& & X \ar[rr]_-{g} & & Y }$$
and recall that $\nu_g = g' \circ \nu_{\id_Y}.$ By the universal property of $\cN(\id_Y),$ we get
a unique arrow $\overline{g}$ such that $\overline{g} \cdot n_{\id_Y} = n_{\id_X} \cdot g$ and 
$\overline{g} \circ \nu_{\id_Y} = \nu_{\id_X} \circ g.$ Because of the first condition on $\overline{g},$ we can use the
universal property of the pullback $\cN(g)$ and we get a unique arrow $v$ such that $v \cdot g' = \overline{g}$ and 
$v \cdot n_g = n_{\id_X}.$ Since
$$v \circ \nu_g = v \cdot g' \circ \nu_{\id_Y} = \overline{g} \circ \nu_{\id_Y} = \nu_{\id_X} \circ g$$
the fact that the $\Theta$-kernel $\cN(g)$ is strong gives us a unique nullhomotopy $\tau \in \Theta(v)$ such that 
$\tau \circ n_g = \nu_{\id_X}.$ Consider now the following diagram, where the solid part is commutative and where
$\overline{z}$ (respectively, $\overline{x}$) is the factorization obtained as in Definition \ref{DefSPropSGlob}.1 when 
we start from the object $\cN(g)$ (respectively, $X$), as in the second part of the proof of Proposition \ref{PropExactSnail1}
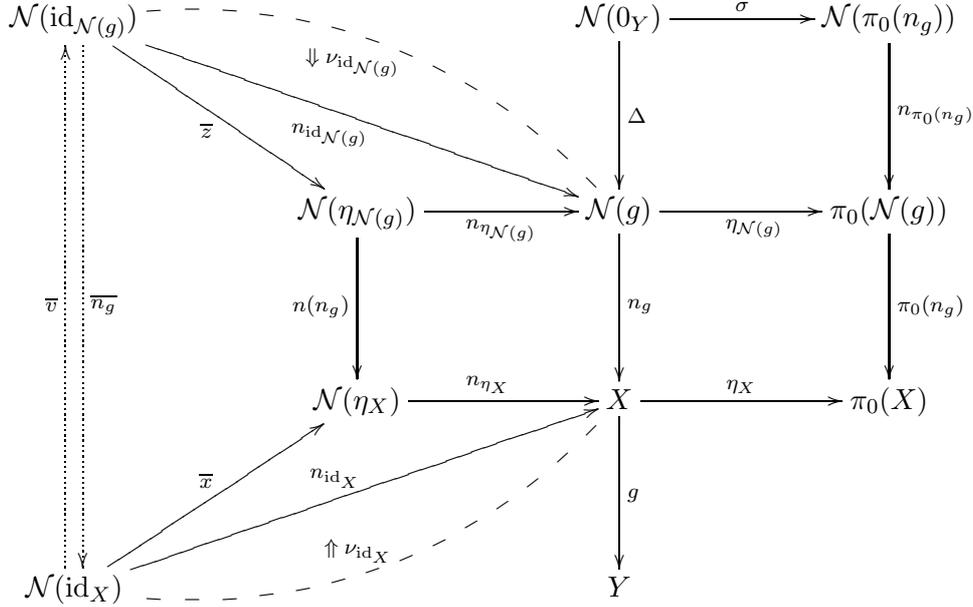
$$\xymatrix{\cN(\id_{\cN(g)}) \ar[rrdd]_{\overline{z}} \ar[rrrrdd]_{n_{\id_{\cN(g)}}} \ar@{.>}@<0.7ex>[dddddd]^{\overline{n_g}} 
\ar@{--}@/^2.5pc/[rrrrdd]_{\Downarrow \; \nu_{\id_{\cN(g)}}}
& & & & \cN(0_Y) \ar[dd]^{\Delta} \ar[rr] ^-{\sigma} & & \cN(\pi_0(n_g)) \ar[dd]^{n_{\pi_0(n_g)}} \\ \\
& & \cN(\eta_{\cN(g)}) \ar[rr]_-{n_{\eta_{\cN(g)}}} \ar[dd]_{n(n_g)} & & \cN(g) \ar[dd]^{n_g} \ar[rr]_{\eta_{\cN(g)}} 
& & \pi_0(\cN(g)) \ar[dd]^{\pi_0(n_g)} \\ \\
& & \cN(\eta_X) \ar[rr]^{n_{\eta_X}} & & X \ar[rr]^-{\eta_X} \ar[dd]^{g} & & \pi_0(X) \\ \\
\cN(\id_X) \ar@{.>}@<0.7ex>[uuuuuu]^{\overline{v}} \ar[rruu]^{\overline{x}} \ar[rrrruu]^{n_{\id_X}} 
\ar@{--}@/_2.5pc/[rrrruu]^{\Uparrow \; \nu_{\id_X}} & & & & Y }$$
By the universal property of $\cN(\id_{\cN(g)}),$ we get a unique arrow $\overline{v}$ such that
$\overline{v} \cdot n_{\id_{\cN(g)}} = v$ and $\overline{v} \circ \nu_{\id_{\cN(g)}} = \tau.$ Moreover, the universal property
of $\cN(\id_X)$ gives a unique arrow $\overline{n_g}$ such that $\overline{n_g} \cdot n_{\id_X} = n_{\id_{\cN(g)}} \cdot n_g$
and $\overline{n_g} \circ \nu_{\id_X} = \nu_{\id_{\cN(g)}} \circ n_g.$ Now we check thet $\overline{v}$ and $\overline{n_g}$
realize an isomorphism. We will use three times the cancellation property recalled in \ref{TextCanc}. First, we check that 
$\overline{n_g} \cdot v = n_{\id_{\cN(g)}} \colon$
\begin{enumerate}
\item[-] $\overline{n_g} \cdot v \cdot n_g = \overline{n_g} \cdot n_{\id_X} = n_{\id_{\cN(g)}} \cdot n_g$
\item[-] $\overline{n_g} \cdot v \circ \nu_g = \overline{n_g} \circ \nu_{\id_X} \circ g = \nu_{\id_{\cN(g)}} \circ n_g \cdot g = 
n_{\id_{\cN(g)}} \circ \nu_g,$ the last equality coming from the reduced interchange
\end{enumerate}
Second, we check that $\overline{n_g} \cdot \overline{v} = \id \colon$
\begin{enumerate}
\item[-] $\overline{n_g} \cdot \overline{v} \cdot n_{\id_{\cN(g)}} = \overline{n_g} \cdot v = n_{\id_{\cN(g)}}$
\item[-] $\overline{n_g} \cdot \overline{v} \circ \nu_{\id_{\cN(g)}} = \overline{n_g} \circ \tau = \nu_{\id_{\cN(g)}}$ where, for the last equality, we 
compose with $n_g \colon$
\item[-] $\overline{n_g} \circ \tau \circ n_g = \overline{n_g} \circ \nu_{\id_X} = \nu_{\id_{\cN(g)}} \circ n_g$
\end{enumerate}
Third, we check that $\overline{v} \cdot \overline{n_g} = \id \colon$
\begin{enumerate}
\item[-] $\overline{v} \cdot \overline{n_g} \cdot n_{\id_X} = \overline{v} \cdot n_{\id_{\cN(g)}} \cdot n_g = v \cdot n_g = n_{\id_X}$
\item[-] $\overline{v} \cdot \overline{n_g} \circ \nu_{\id_X} = \overline{v} \circ \nu_{\id_{\cN(g)}} \circ n_g = \tau \circ n_g = \nu_{\id_X}$
\end{enumerate}
It remains to check that $\overline{z} \cdot n(n_g) = \overline{n_g} \cdot \overline{x}.$ For this, it is enough to compose with 
$n_{\eta_X}$ which is a monomorphism because $\pi_0(X)$ is discrete:
\begin{enumerate}
\item[-]  $\overline{z} \cdot n(n_g) \cdot n_{\eta_X} = \overline{z} \cdot n_{\eta_{\cN(g)}} \cdot n_g = n_{\id_{\cN(g)}} \cdot n_g =
\overline{n_g} \cdot n_{\id_X} = \overline{n_g} \cdot \overline{x} \cdot n_{\eta_X}$
\end{enumerate} 
We can conclude as follows: since $X$ is $\cS$-proper and $\overline{n_g}$ is an isomorphism, the equality
$\overline{z} \cdot n(n_g) = \overline{n_g} \cdot \overline{x}$ implies that $n(n_g)$ is in $\cS$ by Condition \ref{CondSurjLike}. 
Since $\cN(g)$ is $\cS$-global and $\Delta \colon \cN(0_Y) \to \cN(g)$ is the kernel of $n_g$ (see Lemma \ref{LemmaDelta}), 
we can use Condition \ref{CondB} and $\sigma$ is in $\cS.$
\end{proof}

\begin{Example}\label{ExArr4}{\rm
Let us consider once again $(\cB,\Theta) = (\Arr(\cA),\Theta_{\Delta})$ and take as $\cS$ the class of arrows
$$\xymatrix{X \ar[rr]^-{g} \ar[d]_{x} & & Y \ar[d]^{y} \\ X_0 \ar[rr]_{g_0} & & Y_0}$$
such that both $g$ and $g_0$ are regular epimorphisms in $\cA.$
To start, assume just that $\cA$ has a zero object, kernels and cokernels. Under these conditions, $\cS$ contains isomorphisms,
each object is $\cS$-global and an object $(Y,y,Y_0)$ is $\cS$-proper precisely when the factorization of $y$ through the kernel 
of its cokernel is a regular epimorphism. This is the definition of proper arrow used in \cite{DB3x3, EVSnail}. Now keep in
mind that limits and colimits in $\Arr(\cA)$ are computed level-wise in $\cA.$ If we add the assumption that $\cA$ is regular,
then $\cS$ satisfies Condition \ref{CondSurjLike}. Finally, if we assume that $\cA$ is regular and protomodular, then $\cS$ 
satisfies also Condition \ref{CondB}, as proved in \cite{DB3x3}. In particular, if $\cA$ is abelian, then $\cS$ satisfies Conditions
\ref{CondSurjLike} and \ref{CondB} and each object in $\Arr(\cA)$ is $\cS$-global and $\cS$-proper. In conclusion, the exactness
of the snail sequence appearing in \cite{ZJEVSnail,EVSnail} (see Example \ref{ExArr2}) is the special case of Propositions 
\ref{PropExactSnail1} and \ref{PropExactSnail2} when $(\cB,\Theta) = (\Arr(\cA),\Theta_{\Delta})$ and $\cA$ is pointed, regular 
and protomodular.
}\end{Example} 

\begin{Remark}\label{RemProofSimilar}{\rm
The proof of the exactness of the snail sequence (Propositions \ref{PropExactSnail1} and \ref{PropExactSnail2}) is more elaborate 
but essentially very similar to the proof of the exactness of the snail sequence in the context of pointed regular protomodular 
categories done in \cite{EVSnail}. Strange enough, Lemma \ref{LemmaDelta} , which is essential in the present proof, does not 
appear in \cite{EVSnail} but comes from the snail lemma for internal groupoids established in \cite{MMVSnail}. This suggests that 
something is still to be understood concerning the generality of the homotopical version of the snail lemma.
}\end{Remark} 

\section{Sequentiable families of arrows}\label{SecSeqFam}

In this section, we fix a category $\cA$ with a zero object 0, kernels and cokernels. We start with the definition of sequentiable 
families of arrows in $\cA,$ with their morphisms and nullhomotopies.

\begin{Definition}\label{DefSeqFam}{\rm
A sequentiable family of arrows $h_{\bullet}$ is a family of pairs of arrows 
$$h_{\bullet} = \{h_n,i_n \}_{n \in \mathbb Z}$$ 
in $\cA$ with $i_n$ connecting the cokernel of $h_{n+1}$ with the kernel of $h_n \colon$
$$\ldots \xymatrix{\Cod(h_{n+1}) \ar[r]^-{q_{n+1}} & \Cok(h_{n+1}) \ar[r]^-{i_n} & \Ker(h_n) \ar[r]^-{k_n} & \Dom(h_n) \ar[r]^-{h_n} & \Cod(h_n)} \ldots$$
A morphism of sequentiable families $f_{\bullet} \colon h_{\bullet} \to h'_{\bullet}$ is a family of pairs of arrows
$$f_{\bullet} = \{ \overline{f}_n \colon \Dom(h_n) \to \Dom(h'_n), \underline{f}_n \colon \Cod(h_n) \to \Cod(h'_n)\}_{n \in \mathbb Z}$$
such that for all $n \in \mathbb Z$
$$\overline{f}_n \cdot h'_n = h_n \cdot \underline{f}_n \;\;\mbox{ and }\;\; 
\underline{f}_{n+1} \cdot q'_{n+1} \cdot i'_n \cdot k'_n = q_{n+1} \cdot i_n \cdot k_n \cdot \overline{f}_n$$
or, equivalently, such that
$$\overline{f}_n \cdot h'_n = h_n \cdot \underline{f}_n \;\;\mbox{ and }\;\; 
C(f)_{n+1} \cdot i'_n = i_n \cdot K(f)_n$$
where $C(f)_{n+1}$ is the unique arrow such that $\underline{f}_{n+1} \cdot q'_{n+1} = q_{n+1} \cdot C(f)_{n+1}$ and $K(f)_n$
is the unique arrow such that $K(f)_n \cdot k'_n = k_n \cdot \overline{f}_n,$ see the next diagram
$$\xymatrix{{} \ar@{}[d]|{\vdots} & & {} \ar@{}[d]|{\vdots} \\ 
\Dom(h_{n+1}) \ar[rr]^-{\overline{f}_{n+1}} \ar[d]_{h_{n+1}} & & \Dom(h'_{n+1}) \ar[d]^{h'_{n+1}} \\
\Cod(h_{n+1}) \ar[rr]^-{\underline{f}_{n+1}} \ar[d]_{q_{n+1}} & & \Cod(h'_{n+1}) \ar[d]^{q'_{n+1}} \\
\Cok(h_{n+1}) \ar@{.>}[rr]^-{C(f)_{n+1}} \ar[d]_{i_n} & & \Cok(h'_{n+1}) \ar[d]^{i'_n} \\
\Ker(h_n) \ar@{.>}[rr]^-{K(f)_n} \ar[d]_{k_n} & & \Ker(h'_n) \ar[d]^{k'_n} \\
\Dom(h_n) \ar[rr]^-{\overline{f}_n} \ar[d]_{h_n} & & \Dom(h'_n) \ar[d]^{h'_n} \\
\Cod(h_n) \ar[rr]^-{\underline{f}_n} & & \Cod(h'_n) \\
{} \ar@{}[u]|{\vdots} & & {} \ar@{}[u]|{\vdots} }$$
With the obvious identity morphisms and composition of morphisms, sequentiable families and their morphisms give rise to a 
category denoted $\Seq(\cA).$
}\end{Definition} 

\begin{Remark}\label{RemKerSeq}{\rm
Observe that the category $\Seq(\cA)$ has kernels and cokernels constructed level-wise in $\cA.$ The connecting arrows are 
obtained by an easy diagram chasing left to the reader. More in general, $\Seq(\cA)$ inherits level-wise all limits and colimits
which eventually exist in $\cA.$
}\end{Remark} 

The category $\Seq(\cA)$ of sequentiable families is equipped with a structure of nullhomotopies which extends the one we 
considered in $\Arr(\cA),$ see Example \ref{ExArr1}.

\begin{Definition}\label{DefNullSeq}{\rm
Let $f_{\bullet} \colon h_{\bullet} \to h'_{\bullet}$ be a morphism of sequentiable families. A nullhomotopy 
$\lambda_{\bullet} \in \Theta_{\Delta}(f_{\bullet})$ is a family of arrows 
$$\lambda_{\bullet} = \{ \lambda_n \colon \Cod(h_n) \to \Dom(h'_n) \}_{n \in \mathbb Z}$$
such that $h_n \cdot \lambda_n = \overline{f}_n$ and $\lambda_n \cdot h'_n = \underline{f}_n$ for all $n \in \mathbb Z$
$$\xymatrix{\Dom(h_n) \ar[rr]^-{\overline{f}_n} \ar[d]_{h_n} & & \Dom(h'_n) \ar[d]^{h'_n} \\
\Cod(h_n) \ar[rr]_-{\underline{f}_n} \ar[rru]^{\lambda_n} & & \Cod(h'_n) }$$
The composition of a nullhomotopy with morphisms on the left and on the right is defined level-wise as in $(\Arr(\cA),\Theta_{\Delta}).$
}\end{Definition} 

\begin{Proposition}\label{PropNullSeq}
Consider the category with nullhomotopies $(\Seq(\cA),\Theta_{\Delta})$ as in Definitions \ref{DefSeqFam} and \ref{DefNullSeq}.
\begin{enumerate}
\item For a morphism $f_{\bullet} \in \Seq(\cA),$ if $\Theta_{\Delta}(f_{\bullet}) \neq \emptyset$ then for all $n \in \mathbb Z$ the arrows
$K(f)_n$ and $C(f)_n$ are zero arrows.
\item The structure $\Theta_{\Delta}$ satisfies the reduced interchange condition \ref{CondRedInter}. 
\item The sequentiable family $0_{\bullet} = \{ \id_0, \id_0 \}_{n \in \mathbb Z}$ is a $\Theta_{\Delta}$-strong zero object of $\Seq(\cA).$
\item If $\cA$ has pullbacks, then $\Seq(\cA)$ has $\Theta_{\Delta}$-kernels constructed level-wise as in $\Arr(\cA),$ and they are strong.
\item Dually, if $\cA$ has pushouts, then $\Seq(\cA)$ has $\Theta_{\Delta}$-cokernels constructed level-wise as in $\Arr(\cA),$ and they are strong.
\end{enumerate}
\end{Proposition}

\begin{proof}
The only thing that deserves some comments is the construction of the connecting arrows in the $\Theta_{\Delta}$-kernel of a morphism.
The $\Theta_{\Delta}$-kernel of a morphism $f_{\bullet} \colon h_{\bullet} \to h'_{\bullet}$ in $\Seq(\cA)$ being constructed level-wise as in 
$\Arr(\cA)$ (see Example \ref{ExArr1}), we adopt the notation depicted in the following diagram (where the region marked as p.b. is a pullback):
$$\xymatrix{\Dom(h_n) \ar[rrrr]^-{\overline{f}_n} \ar[rrd]^{h^P_n} \ar[dd]_{h_n} & & & & \Dom(h'_n) \ar[dd]^{h'_n} \\
& & P_n \ar[rru]_{\pi'_n} \ar[lld]^{\pi_n} \ar@{}[rrd]|{\mbox{p.b.}} \\
\Cod(h_n) \ar[rrrr]_-{\underline{f}_n} & & & & \Cod(h'_n) }$$
Consider now the following commutative diagram, which is the $\Theta_{\Delta}$-kernel $(\cN(f_{\bullet}),n_{f_{\bullet}},\nu_{f_{\bullet}})$
of $f_{\bullet}$ from level $n$+1 to level $n \colon$
$$\xymatrix{\Dom(h_{n+1}) \ar[rr]^-{\id} \ar[d]_{h^P_{n+1}} & & \Dom(h_{n+1}) \ar[rr]^-{\overline{f}_{n+1}} \ar[d]^>>>{h_{n+1}}  
& & \Dom(h'_{n+1}) \ar[d]^{h'_{n+1}} \\
P_{n+1} \ar[rr]_-{\pi_{n+1}} \ar[d]_{q^P_{n+1}} \ar@{-->}[rrrru]^<<<<<<<<<<<<<<<<{\pi'_{n+1}} 
& & \Cod(h_{n+1}) \ar[rr]_-{\underline{f}_{n+1}} \ar[d]^{q_{n+1}} & & \Cod(h'_{n+1}) \ar[d]^{q'_{n+1}} \\
\Cok(h^P_{n+1}) \ar[rr]^-{C(\pi)_{n+1}} \ar@{.>}[d]_{i^P_n} & & \Cok(h_{n+1}) \ar[rr]^-{C(f)_{n+1}} \ar[d]^{i_n} & & \Cok(h'_{n+1}) \ar[d]^{i'_n} \\
\Ker(h^P_n) \ar[rr]_-{K(\pi)_n} \ar[d]_{k^P_n} & & \Ker(h_n) \ar[rr]_-{K(f)_n} \ar[d]^{k_n} & & \Ker(h'_n) \ar[d]^{k'_n} \\
\Dom(h_n) \ar[rr]^-{\id} \ar[d]_{h^P_n} & & \Dom(h_n) \ar[rr]^-{\overline{f}_n} \ar[d]^>>>{h_n} & & \Dom(h'_n) \ar[d]^{h'_n} \\
P_n \ar[rr]_-{\pi_n} \ar@{-->}[rrrru]^<<<<<<<<<<<<<<<<<{\pi'_n} & & \Cod(h_n) \ar[rr]_-{\underline{f}_n} & & \Cod(h'_n) }$$
We are going to prove that, for all $n \in \mathbb Z,$ there exists a unique arrow $i^P_n$ such that $i^P_n \cdot K(\pi)_n = C(\pi)_{n+1} \cdot i_n.$ 
This makes $\cN(f_{\bullet}) = \{ h^P_n, i^P_n \}_{n \in \mathbb Z}$ an object of $\Seq(\cA).$ \\
Existence: first, we check that $C(\pi)_{n+1} \cdot i_n \cdot k_n \cdot h^P_n = 0.$ For this, we compose with the projections of the pullback $P_n \colon$
\begin{enumerate}
\item[-] $C(\pi)_{n+1} \cdot i_n \cdot k_n \cdot h^P_n \cdot \pi_n = C(\pi)_{n+1} \cdot i_n \cdot k_n \cdot h_n = C(\pi)_{n+1} \cdot i_n \cdot 0 = 0$
\item[-] $C(\pi)_{n+1} \cdot i_n \cdot k_n \cdot h^P_n \cdot \pi'_n = C(\pi)_{n+1} \cdot i_n \cdot k_n \cdot \overline{f}_n =
 C(\pi)_{n+1} \cdot i_n \cdot K(f)_n \cdot k'_n = C(\pi)_{n+1} \cdot C(f)_{n+1} \cdot i'_n \cdot k'_n = 0$ where, to justify the last equality, we precompose
 with the epimorphism $q^P_{n+1} \colon$
 \item[-] $ q^P_{n+1} \cdot C(\pi)_{n+1} \cdot C(f)_{n+1} \cdot i'_n \cdot k'_n = \pi_{n+1} \cdot \underline{f}_{n+1} \cdot q'_{n+1} \cdot i'_n \cdot k'_n =
 \pi'_{n+1} \cdot h'_{n+1} \cdot q'_{n+1} \cdot i'_n \cdot k'_n = \pi'_{n+1} \cdot 0 \cdot i'_n \cdot k'_n = 0$
\end{enumerate}
Now, from $C(\pi)_{n+1} \cdot i_n \cdot k_n \cdot h^P_n = 0,$ we get a unique arrow $i^P_n$ such that 
$i^P_n \cdot k^P_n = C(\pi)_{n+1} \cdot i_n \cdot k_n.$ This can be rewritten as $i^P_n \cdot K(\pi)_n \cdot k_n = C(\pi)_{n+1} \cdot i_n \cdot k_n.$
Since $k_n$ is a monomorphism, we can conclude that $i^P_n \cdot K(\pi)_n = C(\pi)_{n+1} \cdot i_n.$ \\
Uniqueness: this easily follows composing with the monomorphism $k^P_n = K(\pi)_n \cdot k_n.$
\end{proof} 

\begin{Text}\label{TextPartHK}{\rm
Based on the description of the $\Theta_{\Delta}$-kernel in $\Seq(\cA)$ contained in the proof of Proposition \ref{PropNullSeq}, 
we list here some special cases of $\Theta_{\Delta}$-kernels involved in the snail sequence.
\begin{enumerate}
\item The $n$-th level of the $\Theta_{\Delta}$-kernel of $\id_{h_{\bullet}}$ is given by
$$\xymatrix{\Dom(h_n) \ar[d]_{\id} \ar[rr]^-{\id} & & \Dom(h_n) \ar[d]^>>>{h_n} \ar[rr]^-{\id} & & \Dom(h_n) \ar[d]^{h_n} \\
\Dom(h_n) \ar[rr]_-{h_n} \ar@{-->}[rrrru]^<<<<<<<<<<<{\id} & & \Cod(h_n) \ar[rr]_-{\id} & & \Cod(h_n) }$$
and the $n$-th connecting arrow of $\cN(\id_{h_{\bullet}})$ is $\id \colon 0 \to 0.$
\item The $n$-th level of the $\Theta_{\Delta}$-kernel of $0_{h_{\bullet}}$ is given by
$$\xymatrix{0 \ar[d] \ar[rr] & & 0 \ar[d] \ar[rr] & & \Dom(h_n) \ar[d]^{h_n} \\
\Ker(h_n) \ar[rr] \ar@{-->}[rrrru]^<<<<<<<<<{k_{h_n}} & & 0 \ar[rr] & & \Cod(h_n) }$$
and the $n$-th connecting arrow of $\cN(0_{h_{\bullet}})$ is $0 \colon \Ker(h_{n+1}) \to 0.$
\item The $n$-th level of the morphism $\eta_{h_{\bullet}} \colon h_{\bullet} \to \pi_0(h_{\bullet})$ is given by
$$\xymatrix{\Dom(h_n) \ar[rr] \ar[d]_{h_n} & & 0 \ar[d] \\
\Cod(h_n) \ar[rr]_-{c_{h_n}} & & \Cok(h_n) }$$
and the $n$-th connecting arrow of $\pi_0(h_{\bullet})$ is $0 \colon \Cok(h_{n+1}) \to 0.$
\item The $n$-th level of the factorization $\overline{h}_{\bullet} \colon \cN(\id_{h_{\bullet}}) \to \cN(\eta_{h_{\bullet}})$
of $n_{\id_{h_{\bullet}}} \colon \cN(\id_{h_{\bullet}}) \to h_{\bullet}$ through 
$n_{\eta_{h_{\bullet}}} \colon \cN(\eta_{h_{\bullet}}) \to h_{\bullet}$ (see Definition \ref{DefSPropSGlob}) is given by
$$\xymatrix{\Dom(h_n) \ar[rr]^-{\id} \ar[d]_{\id} & & \Dom(h_n) \ar[d]^{\overline{h}_n} \\
\Dom(h_n) \ar[rr]_-{\overline{h}_n} & & \Ker(c_{h_n}) }$$
where $\overline{h}_n$ is the factorization of $h_n$ through the kernel of its cokernel.
\end{enumerate}
}\end{Text} 

\begin{Text}\label{TextSnailSeq}{\rm
Assume now that the category $\cA$ has pullbacks. 
Thanks to Proposition \ref{PropNullSeq}, we can perform the construction of Section \ref{SecHomSnailSeq} in the category 
with nullhomotopies $(\Seq(\cA),\Theta_{\Delta}).$ Starting from a morphism $f_{\bullet} \colon h_{\bullet} \to h'_{\bullet},$ 
we get a sequence
$$\xymatrix{\cN(0_{\cN(f_{\bullet})}) \ar[r] & \cN(0_{h_{\bullet}}) \ar[r] & \cN(0_{h'_{\bullet}}) \ar[r] & 
\pi_0(\cN(f_{\bullet})) \ar[r] & \pi_0(h_{\bullet}) \ar[r] & \pi_0(h'_{\bullet}) }$$
Keeping in mind the various special cases of $\Theta_{\Delta}$-kernels described in \ref{TextPartHK}, we can make the 
previous sequence explicit. Here it is from level $n$+1 to level $n \colon$
$$\resizebox{\displaywidth}{!}{
\xymatrix{0 \ar[r] \ar[d] & 0 \ar[r] \ar[d] & 0 \ar[rr] \ar[d] & & 0 \ar[r] \ar[d] & 0 \ar[r] \ar[d] & 0 \ar[d] \\
\Ker(h^P_{n+1}) \ar[r]^-{K(\pi)_{n+1}} \ar[d]_{\id} & \Ker(h_{n+1}) \ar[r]^-{K(f)_{n+1}} \ar[d]_{\id} & 
\Ker(h'_{n+1}) \ar[rr]^-{\langle 0,k'_{n+1} \rangle \cdot q'_{n+1}} \ar[d]_{\id} & &
 \Cok(h^P_{n+1}) \ar[r]^-{C(\pi)_{n+1}} \ar[d]^{\id} & \Cok(h_{n+1}) \ar[r]^-{C(f)_{n+1}} \ar[d]^{\id} & \Cok(h'_{n+1}) \ar[d]^{\id} \\
\Ker(h^P_{n+1}) \ar[r]^-{K(\pi)_{n+1}} \ar[d] & \Ker(h_{n+1}) \ar[r]^-{K(f)_{n+1}} \ar[d] & 
\Ker(h'_{n+1}) \ar[rr]^-{\langle 0,k'_{n+1} \rangle \cdot q'_{n+1}} \ar[d] & & 
\Cok(h^P_{n+1}) \ar[r]^-{C(\pi)_{n+1}} \ar[d] & \Cok(h_{n+1}) \ar[r]^-{C(f)_{n+1}} \ar[d] & \Cok(h'_{n+1}) \ar[d] \\
0 \ar[r] \ar[d] & 0 \ar[r] \ar[d] & 0 \ar[rr] \ar[d] & & 0 \ar[r] \ar[d] & 0 \ar[r] \ar[d] & 0 \ar[d] \\
0 \ar[r] \ar[d] & 0 \ar[r] \ar[d] & 0 \ar[rr] \ar[d] & & 0 \ar[r] \ar[d] & 0 \ar[r] \ar[d] & 0 \ar[d] \\
\Ker(h^P_n) \ar[r]^-{K(\pi)_n} & \Ker(h_n) \ar[r]^-{K(f)_n} & \Ker(h'_n) \ar[rr]^-{\langle 0,k'_{n} \rangle \cdot q'_{n}} 
& & \Cok(h^P_n) \ar[r]^-{C(\pi)_n} & \Cok(h_n) \ar[r]^-{C(f)_n} & \Cok(h'_n) }}$$
}\end{Text} 

\begin{Text}\label{TextExExn}{\rm
Consider now a class $\cS$ of arrows in $\cA$ and extend it to a classe, still denoted by $\cS,$ of morphisms in 
$\Seq(\cA)$ taking as $\cS$-morphisms in $\Seq(\cA)$ those morphisms $f_{\bullet} \colon h_{\bullet} \to h'_{\bullet}$ 
such that both $\overline{f}_n$ and $\underline{f}_n$ are in $\cS$ for all $n \in \mathbb Z.$ Since all the objects involved 
are discrete, and thanks to Remark \ref{RemKerSeq}, the $\cS$-exactness in $\Seq(\cA)$ of the sequence in \ref{TextSnailSeq}
amounts to the $\cS$-exactness in $\cA,$ for all $n \in \mathbb Z,$  of the sequence
$$\xymatrix{\Ker(h^P_n) \ar[r]^-{K(\pi)_n} & \Ker(h_n) \ar[r]^-{K(f)_n} & \Ker(h'_n) \ar[rr]^-{\langle 0,k'_{n} \rangle \cdot q'_{n}} 
& & \Cok(h^P_n) \ar[r]^-{C(\pi)_n} & \Cok(h_n) \ar[r]^-{C(f)_n} & \Cok(h'_n) }$$
}\end{Text} 

\begin{Text}\label{TextPasting}{\rm
Let us observe also that the connecting arrows of the sequentiable families $\cN(f_{\bullet}), h_{\bullet}$ and $h'_{\bullet}$ allow us to past
together all the six-term sequences in $\cA,$ getting a long sequence as follows:
$$\resizebox{\displaywidth}{!}{
\xymatrix{\Ker(h^P_{n+1}) \ar[r]^-{K(\pi)_{n+1}} & \Ker(h_{n+1}) \ar[r]^-{K(f)_{n+1}} & \Ker(h'_{n+1}) \ar[rr]^-{\langle 0,k'_{{n+1}} \rangle \cdot q'_{{n+1}}} 
& & \Cok(h^P_{n+1}) \ar[r]^-{C(\pi)_{n+1}} \ar[lllldd] _{i^P_n} & \Cok(h_{n+1}) \ar[r]^-{C(f)_{n+1}} \ar[lllldd]_{i_n} & \Cok(h'_{n+1}) \ar[lllldd]_{i'_n} \\ \\ 
\Ker(h^P_n) \ar[r]_-{K(\pi)_n} & \Ker(h_n) \ar[r]_-{K(f)_n} & \Ker(h'_n) \ar[rr]_-{\langle 0,k'_{n} \rangle \cdot q'_{n}} 
& & \Cok(h^P_n) \ar[r]_-{C(\pi)_n} & \Cok(h_n) \ar[r]_-{C(f)_n} & \Cok(h'_n) }}$$
}\end{Text} 

Here is an ad hoc condition which ensures the preservation of exactness of the six-term sequences in $\cA$ whene we past 
them all together.

\begin{Definition}\label{DefIsoSeq}{\rm
We denote by $\IsoSeq(\cA)$ the full subcategory of $\Seq(\cA)$ of the isosequentiable families, that is, those families
$h_{\bullet} = \{h_n,i_n\}_{n \in \mathbb Z}$ such that $i_n$ is an isomorphism for all $n \in \mathbb Z.$
}\end{Definition} 

\begin{Corollary}\label{CorLongSeq}
Let $f_{\bullet} \colon h_{\bullet} \to h'_{\bullet}$ be a morphism in $\IsoSeq(\cA).$ 
Take as $\cS$ the class of regular epimorphisms in $\cA$ and extend $\cS$ to $\Seq(\cA)$ as in \ref{TextExExn}.
Consider the following sequences, the first one being in $\Seq(\cA)$ and the second one being in $\cA \colon$
$$\xymatrix{\cN(0_{\cN(f_{\bullet})}) \ar[r] & \cN(0_{h_{\bullet}}) \ar[r] & \cN(0_{h'_{\bullet}}) \ar[r] & 
\pi_0(\cN(f_{\bullet})) \ar[r] & \pi_0(h_{\bullet}) \ar[r] & \pi_0(h'_{\bullet}) }$$
$$\resizebox{\displaywidth}{!}{
\xymatrix{\Ker(h^P_{n+1}) \ar[r]^-{K(\pi)_{n+1}} & \Ker(h_{n+1}) \ar[r]^-{K(f)_{n+1}} & \Ker(h'_{n+1}) \ar[rr]^-{\langle 0,k'_{{n+1}} \rangle \cdot q'_{{n+1}}} 
& & \Cok(h^P_{n+1}) \ar[r]^-{C(\pi)_{n+1}} & \Cok(h_{n+1}) \ar[r]^-{C(f)_{n+1}} \ar[llld]|{C(f)_{n+1} \cdot i'_n = i_n \cdot K(f)_n} & \Cok(h'_{n+1}) \\
\Ker(h^P_n) \ar[r]_-{K(\pi)_n} & \Ker(h_n) \ar[r]_-{K(f)_n} & \Ker(h'_n) \ar[rr]_-{\langle 0,k'_{n} \rangle \cdot q'_{n}} 
& & \Cok(h^P_n) \ar[r]^-{C(\pi)_n} & \Cok(h_n) \ar[r]^-{C(f)_n} \ar[llld]|{C(f)_{n} \cdot i'_{n-1} = i_{n-1} \cdot K(f)_{n-1}} & \Cok(h'_n) \\
\Ker(h^P_{n-1}) \ar[r]_-{K(\pi)_{n-1}} & \Ker(h_{n-1}) \ar[r]_-{K(f)_{n-1}} & \Ker(h'_{n-1}) \ar[rr]_-{\langle 0,k'_{{n-1}} \rangle \cdot q'_{{n-1}}} 
& & \Cok(h^P_{n-1}) \ar[r]_-{C(\pi)_{n-1}} & \Cok(h_{n-1}) \ar[r]_-{C(f)_{n-1}} & \Cok(h'_{n-1})}}$$
\begin{enumerate}
\item If $\cA$ is regular and protomodular and if $\cN(f_{\bullet}), h_{\bullet}$ and $h'_{\bullet}$ are $\cS$-proper, 
then the six-term sequence in $\Seq(\cA)$ and each row of the long sequence in $\cA$ are $S$-exact.
\item If moreover $h_{\bullet}$ and $h_{\bullet}'$ are isosequentiable, then the long sequence in $\cA$ is $S$-exact at each point.
\end{enumerate} 
\end{Corollary} 

\begin{proof}
1. By Propositions \ref{PropExactSnail1} and \ref{PropExactSnail2}, the six-term sequence in $\Seq(\cA)$ is $S$-exact.
Therefore, by point \ref{TextExExn}, each horizontal row of the long sequence is $\cS$-exact. 
Note that we can use Propositions \ref{PropExactSnail1} and \ref{PropExactSnail2} because the class $\cS$ satisfies Conditions 
\ref{CondSurjLike} and \ref{CondB}. The argument for $\Seq(\cA)$ is the same as for $\Arr(\cA),$ see Example \ref{ExArr4}. \\
2. By assumption, the connecting arrows $i_n$ and $i'_n$ are now isomorphisms. Therefore, we can paste vertically the different 
rows and we are done. 
\end{proof} 

\begin{Remark}\label{remCorLongSeq}{\rm
As already observed in Example \ref{ExArr4} in the case $(\Arr(\cA),\Theta_{\Delta}),$ each object of $\Seq(\cA)$ is $\cS$-global 
and, if $\cA$ is abelian, it is also $\cS$-proper. Therefore, in the abelian case, the assumptions in Corollary \ref{CorLongSeq} that 
$\cN(h_{\bullet}), h_{\bullet}$ and $h'_{\bullet}$ are $S$-proper are redundant.
}\end{Remark} 

\begin{Remark}\label{RemNonIsoSeq}{\rm
In general, the $\Theta_{\Delta}$-kernel of a morphism between isosequentiable families is not isosequentiable. The second
item of point \ref{TextPartHK} provides a counterexample: the $\Theta_{\Delta}$-kernel of the initial morphism $0 \to h_{\bullet}$
is isosequentaible precisely when each $h_n$ has trivial kernel.
}\end{Remark} 

\section{From chain complexes to sequentiable families}\label{SecFromChToSeqFam}

As in Section \ref{SecSeqFam}, we fix a category $\cA$ with a zero object 0, kernels and cokernels. The category 
of chain complexes in $\cA$ will be denoted by $\Ch(\cA)$ and a typical object will be depicted as
$$C_{\bullet} \colon \;\; \ldots \xymatrix{C_{n+1} \ar[rr]^-{d^C_{n+1}} & & C_n \ar[rr]^-{d^C_n} & & C_{n-1} \ar[rr]^-{d^C_{n-1}} & & C_{n-2} } \ldots$$

\begin{Text}\label{TextPropArr}{\rm
In this section, we fix as class $\cS$ of arrows in $\cA$ the class of regular epimorphisms and we write proper instead of $\cS$-proper.
By Example \ref{ExArr4}, we already know that an arrow in $\cA$ is proper (in the sense that its factorization through the 
kernel of its cokernel is a regular epimorphism) iff it is proper as an object of $(\Arr(\cA),\Theta_{\Delta}).$ Similarly, for a sequentiable family 
$h_{\bullet} = \{h_n,i_n\}_{n \in \mathbb Z},$ to be proper as an object of $(\Seq(\cA),\Theta_{\Delta})$ amounts to the fact that each $h_n$ 
is a proper arrow. By analogy, we will say that a complex $C_{\bullet} \in \Ch(\cA)$ is proper if each $d^C_n$ is a proper 
arrow. In this way, our terminology for complexes agrees with the terminology introduced in \cite{TETVDL}.
}\end{Text} 

\begin{Text}\label{TextInfSeq}{\rm 
Let us start this last section looking more carefully at the classical link between the snake lemma and the long homology sequence in the 
abelian case (see for example \cite{WBook}). Given an extension of complexes
$$\xymatrix{A_{\bullet} \ar[rr]^-{f_{\bullet}} & & B_{\bullet} \ar[rr]^-{g_{\bullet}} & & C_{\bullet} }$$
we get a family of extensions in $\Arr(\cA)$
$$\xymatrix{A_n \ar[rr]^-{f_n} \ar[d]_{d^A_n} & & B_n \ar[rr]^-{g_n} \ar[d]_{d^B_n} & & C_n \ar[d]^{d^C_n} \\
A_{n-1} \ar[rr]_-{f_{n-1}} & & B_{n-1} \ar[rr]_-{g_{n-1}} & & C_{n-1} }$$
that is, $g_n$ is the cokernel of $f_n$ and $f_n$ is the kernel of $g_n,$ and this for each $n \in \mathbb Z.$ Nevertheless, in order to get the 
homology sequence, we do not apply the snake lemma directly to these extensions. Instead, we factorize each one of them 
and we get a new family of dotted commutative diagrams (notation $\cF(-)$ is explained in Proposition \ref{PropSeqAssCh})
$$\resizebox{\displaywidth}{!}{
\xymatrix{A_n \ar[rrrr]^-{f_n} \ar[ddd]_{d^A_n} \ar[rrd]^{q^A_n} & & & & B_n \ar[rrrr]^-{g_n} \ar[ddd]_{d^B_n} \ar[rrd]^{q^B_n} 
& & & & C_n \ar[ddd]_{d^C_n} \ar[rrd]^{q^C_n} \\ 
& & \Cok(d^A_{n+1}) \ar@{.>}[rrrr]^<<<<<<<<<<<<{\overline{\cF f}_n} \ar@{.>}[d]_{h^{\cF A}_n} & & & & 
\Cok(d^B_{n+1}) \ar@{.>}[rrrr]^<<<<<<<<<<<<{\overline{\cF g}_n} \ar@{.>}[d]_{h^{\cF B}_n} & & & & \Cok(d^C_{n+1}) \ar@{.>}[d]_{h^{\cF C}_n} \\ 
& & \Ker(d^A_{n-1}) \ar@{.>}[rrrr]_<<<<<<<<<<<<{\underline{\cF f}_n} \ar[lld]^{k^A_n} & & & & 
\Ker(d^B_{n-1}) \ar@{.>}[rrrr]_<<<<<<<<<<<<{\underline{\cF g}_n} \ar[lld]^{k^B_n} & & & & \Ker(d^C_{n-1}) \ar[lld]^{k^C_n} \\
A_{n-1} \ar[rrrr]_-{f_{n-1}} & & & & B_{n-1} \ar[rrrr]_-{g_{n-1}} & & & & C_{n-1} }}$$
Finally we apply the snake lemma to this second family of diagrams (this is possible because $\overline{\cF g}_n$ is still the 
cokernel of $\overline{\cF f}_n$ and $\underline{\cF f}_n$ is still the kernel of $\underline{\cF g}_n$) and we get the long exact sequence 
in homology. Indeed, the kernel of $h^{\cF A}_n$ is the homology object $H_n(A_{\bullet})$ and the cokernel of $h^{\cF A}_n$ is the 
homology object $H_{n-1}(A_{\bullet}),$ and the same holds for the complexes $B_{\bullet}$ and $C_{\bullet}.$ The sequentiable families 
of arrows in $\cA$ arise precisely from this intermediate construction, as formalized in the next proposition.
}\end{Text} 

\begin{Proposition}\label{PropSeqAssCh}
From any complex $C_{\bullet} \in \Ch(\cA),$ we obtain a sequentiable family of arrows 
$\cF(C_{\bullet}) = \{h^{\cF C}_n, i^{\cF C}_n \}_{n \in \mathbb Z}  \in \Seq(\cA)$ depicted hereunder (from level $n$+1 to level $n$)
$$\xymatrix{C_{n+1} \ar[rrrr]^-{q^C_{n+1}} \ar[rrrrd]_{k(d^C_{n+1})} \ar[ddd]_{d^C_{n+1}} 
& & & & \Cok(d^C_{n+2}) \ar[d]^{h^{\cF C}_{n+1}} \ar[llllddd]_{q(d^C_{n+1})} \\
& & & & \Ker(d^C_n) \ar[d]^{q^{\cF C}_{n+1}} \ar[lllldd]^{k^C_n} \\
& & & & \Cok(h^{\cF C}_{n+1}) \ar[dd]^{i^{\cF C}_n} \\
C_n \ar[ddd]_{d^C_n} \ar[rrrrdd]^{q^C_n} \ar[rrrrddd]_{k(d^C_n)} \\
& & & & \Ker(h^{\cF C}_n) \ar[d]^{k^{\cF C}_n} \\
& & & & \Cok(d^C_{n+1}) \ar[d]^{h^{\cF C}_n} \ar[lllld]_{q(d^C_n)} \\
C_{n-1} & & & & \Ker(d^C_{n-1}) \ar[llll]^-{k^C_{n-1}} }$$
where $k(d^C_n)$ is the unique arrow such that $k(d^C_n) \cdot k^C_{n-1} = d^C_n$ and $q(d^C_n)$ is the unique arrow 
such that $q^C_n \cdot q(d^C_n) = d^C_n.$ The components of the family $\cF(C_{\bullet})$ are as follows:
\begin{enumerate}
\item[-] $h^{\cF C}_n$ is the unique arrow such that $q^C_n \cdot h^{\cF C}_n = k(d^C_n)$ or, equivalently, the unique arrow 
such that $h^{\cF C}_n \cdot k^C_{n-1} = q(d^C_n).$
\item[-] $i^{\cF C}_n$ is the unique arrow such that $q^{\cF C}_{n+1} \cdot i^{\cF C}_n \cdot k^{\cF C}_n = k^C_n \cdot q^C_n.$
\end{enumerate}
This construction extends to a functor
$$\cF \colon \Ch(\cA) \longrightarrow \Seq(\cA)$$
Moreover, if $\cA$ is regular and protomodular, then the family $\cF(C_{\bullet})$ is isosequentiable provided that the complex 
$C_{\bullet}$ is proper. \\
In particular, if $\cA$ is abelian, we get a functor $\cF \colon \Ch(\cA) \to \IsoSeq(\cA).$
\end{Proposition} 

\begin{proof}
The equivalence between the two possible definitions of $h^{\cF C}_n$ comes from the fact that $q^C_n$ is an epimorphism 
and $k^C_{n-1}$ is a monomorphism. As far as $i^{\cF C}_n$ is concerned, observe that, since $q^C_{n+1}$ is an epimorphism, 
$\Cok(h^{\cF C}_{n+1}) = \Cok(k(d^C_{n+1})).$ Similarly, since $k^C_{n-1}$ is a monomoprphism, $\Ker(h^{\cF C}_n) = \Ker(q(d^C_n)).$
Now the argument runs as usual: since $k(d^C_{n+1}) \cdot k^C_n \cdot q^C_n = 0,$ there exists a unique 
$j_n \colon \Cok(k(d^C_{n+1})) \to \Cok(d^C_{n+1})$ such that $q^{\cF C}_{n+1} \cdot j_n = k^C_n \cdot q^C_n.$ Since moreover
$j_n \cdot q(d^C_n) = 0$ (precompose with $q^{\cF C}_{n+1}$), there exists a unique $i^{\cF C}_n$ such that 
$i^{\cF C}_n \cdot k^{\cF C}_n = j_n.$ \\
The construction of a morphism 
$\cF(g_{\bullet}) = \{\overline{\cF g}_n, \underline{\cF g}_n \}_{n \in \mathbb Z} \colon \cF(B_{\bullet}) \to \cF(C_{\bullet})$ from a morphism 
$g_{\bullet} \colon B_{\bullet} \to C_{\bullet}$ in $\Ch(\cA)$ has been already depicted in the last diagram of \ref{TextInfSeq} and the 
functoriality of the construction is obvious. Let us check just that $\cF(g_{\bullet})$ is indeed a morphism in $\Seq(\cA) \colon$
\begin{enumerate}
\item[-] $\underline{\cF g}_{n+1} \cdot q^{\cF C}_{n+1} \cdot i^{\cF C}_n \cdot k^{\cF C}_n =
\underline{\cF g}_{n+1} \cdot k^C_n \cdot q^C_n = k^B_n \cdot g_n \cdot q^C_n =
k^B_n \cdot q^B_n \cdot \overline{\cF g}_n = q_{n+1}^{\cF B} \cdot i^{\cF B}_n \cdot k^{\cF B}_n \cdot \overline{\cF g}_n$
\end{enumerate}
Finally, keeping in mind what we have observed above about $\Cok(h^{\cF C}_{n+1})$ and $\Ker(h^{\cF C}_n),$ the fact that 
$\cF(C_{\bullet})$ is isosequentiable if $C_{\bullet}$ is proper and $\cA$ is regular and protomodular follows from Lemma 4.5.1
in \cite{BB}. The special case follows once again because in an abelian category all arrows are proper. 
\end{proof}

\begin{Text}\label{TextHomObj}{\rm
Let us point out explicitly from the proof of Proposition \ref{PropSeqAssCh} that the homology objects of the complex $C_{\bullet}$ 
can be recovered, for all $n \in \mathbb Z,$ as $H_n(C_{\bullet}) = \Cok(h^{\cF C}_{n+1})$ or, equivalently, as
$H_n(C_{\bullet}) = \Ker(h^{\cF C}_{n}).$ 
In other words, for an isosequentiable family $h_{\bullet},$ we can put
$H_n(h_{\bullet}) = \Cok(h_{n+1})$ or $H_n(h_{\bullet}) = \Ker(h_n).$
In this way, if $\cA$ is regular and protomodular and if we consider proper complexes, the diagram
$$\xymatrix{\Ch(\cA) \ar[rr]^-{\cF} \ar[rd]_{H_n} & & \IsoSeq(\cA) \ar[ld]^{H_n} \\ & \cA }$$
commutes (up to isomorphism) for all $n \in \mathbb Z.$
We can relate this fact to Remark \ref{RemNonIsoSeq}: for a complex $C_{\bullet},$ the $\Theta_{\Delta}$-kernel of the initial morphism
$0 \to \cF(C_{\bullet})$ is isosequentiable precisely when $C_{\bullet}$ is acyclic.
}\end{Text} 

\begin{Lemma}\label{LemmaProPro}
Let $\cA$ be a regular category. 
\begin{enumerate}
\item Consider the following situation
$$\xymatrix{A \ar[r]^-{f} & B \ar[r]^-{g} & C \ar[r]^-{h} & D }$$
If $f \cdot g \cdot h$ is proper, $f$ is an epimorphism and $h$ is a monomorphism, then $g$ is proper. 
\item If a complex $C_{\bullet}$ is proper, then the associated sequentiable family $\cF(C_{\bullet})$ is proper. 
\end{enumerate}
\end{Lemma}

\begin{proof}
1. We split the proof in two parts. \\
(a) We assume that $a = g \cdot h$ is proper and $h$ is a monomorphism and we show that $g$ is proper.
Consider the following commutative diagram
$$\xymatrix{ & & B \ar[rrd]^{\overline{g}} \ar[rrdd]_>>>>>>>>>>{g} \ar[lld]_{\overline{a}} \ar[lldd]^>>>>>>>>>>{a} \\
\Ker(c_a) \ar[d]_{k_{c_a}} \ar@{.>}@<0.5ex>[rrrr]^-{s} & & & & \Ker(c_g) \ar[d]^{k_{c_g}} \ar@{.>}@<0.5ex>[llll]^-{t} \\
D \ar[d]_{c_a} & & & & C \ar[llll]_-{h} \ar[d]^{c_g} \\
\Cok(a) & & & & \Cok(g) \ar[llll]^-{h'} }$$
Since $k_{c_g} \cdot h \cdot c_a = 0,$ there exists a unique arrow $t$ such that $t \cdot k_{c_a} = k_{c_g} \cdot h.$
Since $\overline{a} \cdot k_{c_a} = \overline{g} \cdot k_{c_g} \cdot h,$ with $\overline{a}$ a regular epimorphism
(because $a$ is proper) and $k_{c_g} \cdot h$ a monomorphism, there exists a unique arrow $s$ such that 
$\overline{a} \cdot s = \overline{g}$ and $s \cdot k_{c_g} \cdot h = k_{c_a}.$ Composing on the left with $\overline{a}$
and on the right with $k_{c_a},$ one checks that $s \cdot t = \id.$ Composing on the right with $k_{c_g} \cdot h,$ one checks 
that $t \cdot s = \id.$ Therefore, $\overline{a}$ and $\overline{g}$ are equal up to an isomorphism, so that $\overline{g}$
is a regular epimorphism, that is, $g$ is proper. \\
(b) We assume that $b = f \cdot g$ is proper and $f$ is an epimorphism and we show that $g$ is proper.
Consider the following commutative diagram
$$\xymatrix{A \ar[rrrr]^-{f} \ar[rrdd]^<<<<<<<<<<<<{b} \ar[d]_{\overline{b}} & & & & 
B \ar[d]^{\overline{g}} \ar[lldd]_<<<<<<<<<<<<{g} \\
\Ker(c_b) \ar[rrd]_{k_{c_b}} \ar@{.>}[rrrr]_-{t} & & & & \Ker(c_g) \ar[lld]^{k_{c_g}} \\
& & C \ar[ld]_{c_b} \ar[rd]^{c_g} \\
& \Cok(b) \ar@{.>}[rr]_-{s} & & \Cok(g) }$$
Since $f$ is an epimorphism, the unique arrow $s$ such that $c_b \cdot s = c_g$ is an isomorphism. Therefore, the unique 
arrow $t$ such that $t \cdot k_{c_g} = k_{c_b}$ also is an isomorphism. Since $\overline{b} \cdot t = f \cdot \overline{g}$
(compose with $k_{c_g}$ to check this) and since $\overline{b}$ is a regular epimorphism (because $b$ is proper), we 
conclude that $\overline{g}$ is a regular epimorphism, that is, $g$ is proper. \\
2. This follows from 1. because, for all $n \in \mathbb Z,$ we have $q^C_n \cdot h^{\cF C}_n \cdot k^C_{n-1} = d^C_n,$ 
which is proper (notation as in Proposition \ref{PropSeqAssCh}).
\end{proof}

\begin{Corollary}\label{CorLongExSeq}
Let $\cA$ be regular and protomodular and let $g_{\bullet} \colon B_{\bullet} \to C_{\bullet}$ be a morphism of proper complexes.
If $\cN(\cF(g_{\bullet}))$ is proper, then we get a long exact sequence in homology
$$\ldots \xymatrix{H_{n+1}(B_{\bullet}) \ar[r] & H_{n+1}(C_{\bullet}) \ar[r] & \Cok(h^P_{n+1}) \ar[r] & H_n(B_{\bullet}) \ar[r] & H_n(C_{\bullet}) } \ldots$$
In particular, if $\cA$ is abelian this holds with no assumption on $B_{\bullet}, C_{\bullet}$ and $\cN(\cF(g_{\bullet})).$
\end{Corollary}

\begin{proof}
This is the long exact sequence of Corollary \ref{CorLongSeq} applied to $\cF(g_{\bullet}) \colon \cF(B_{\bullet}) \to \cF(C_{\bullet})$
(we need Lemma \ref{LemmaProPro} to use Corollary \ref{CorLongSeq}).
Therefore, $\{h^P_n\}_{n \in \mathbb Z}$ is the first component of the $\Theta_{\Delta}$-kernel $\cN(\cF(g_{\bullet}))$ in $\Seq(\cA)$
and the homology objects of $B_{\bullet}$ and $C_{\bullet}$ appear in the sequence as explained in \ref{TextHomObj}.
\end{proof}

\begin{Text}\label{TextCompClassSeq}{\rm
It is time to compare the long sequence of Corollary \ref{CorLongExSeq} with the usual one obtained from an extension of complexes,
but this is an easy task. Assume that $\cA$ is regular and protomodular and consider an extension of proper complexes
$$\xymatrix{A_{\bullet} \ar[r]^-{f_{\bullet}} & B_{\bullet} \ar[r]^-{g_{\bullet}} & C_{\bullet} }$$
Assume also that the sequentiable family $\cN(\cF(g_{\bullet}))$ is proper. The long exact sequence of Corollary \ref{CorLongExSeq} is obtained 
pasting together infinitely many copies of the six-term snail sequence coming from 
$$\xymatrix{\Cok(d^B_{n+1}) \ar[rr]^-{\id} \ar[d]_{h^P_n} & & \Cok(d^B_{n+1}) \ar[rr]^-{\overline{\cF g}_n} \ar[d]^{h^{\cF B}_n} 
& & \Cok(d^C_{n+1}) \ar[d]^{h^{\cF C}_n} \\
P_n \ar[rr]_-{\pi_n} & & \Ker(d^B_{n-1}) \ar[rr]_-{\underline{\cF g}_n} & & \Ker(d^C_{n-1}) }$$
The usual long exact sequence in homology is obtained pasting together infinitely many copies of the six-term snake sequence coming from 
$$\xymatrix{\Cok(d^A_{n+1}) \ar[rr]^-{\overline{\cF f}_n} \ar[d]_{h^{\cF A}_n} & & \Cok(d^B_{n+1}) \ar[rr]^-{\overline{\cF g}_n} \ar[d]^{h^{\cF B}_n} 
& & \Cok(d^C_{n+1}) \ar[d]^{h^{\cF C}_n} \\
\Ker(d^A_{n-1}) \ar[rr]_-{\underline{\cF f}_n} & & \Ker(d^B_{n-1}) \ar[rr]_-{\underline{\cF g}_n} & & \Ker(d^C_{n-1}) }$$
Since, for all $n \in \mathbb Z,$ the arrow $\overline{\cF g}_n$ is a regular epimorphism (because 
$g_n \cdot q^C_n = q^B_n \cdot \overline{\cF g}_n,$ see \ref{TextInfSeq}, and $g_n$ is a regular epimorphism), we can apply Proposition 4.3 in 
\cite{EVSnail} and we have that the two six-term exact sequences are isomorphic. We can conclude that, in this particular case, the long exact 
sequence of Corollary \ref{CorLongExSeq} is nothing but the usual one
$$\ldots \xymatrix{H_{n+1}(B_{\bullet}) \ar[r] & H_{n+1}(C_{\bullet}) \ar[r] & H_n(A_{\bullet}) \ar[r] & H_n(B_{\bullet}) \ar[r] & H_n(C_{\bullet}) } \ldots$$
}\end{Text}

\begin{Text}\label{TextCompClassSeqBis}{\rm
Under the same assumptions, we can reformulated point \ref{TextCompClassSeq} in a more synthetic way. First, observe that the functor
$\cF \colon \Ch(\cA) \to \Seq(\cA)$ preserves kernels and regular epimorphisms (use Remark \ref{RemKerSeq}). 
Second, compare the kernel and the $\Theta_{\Delta}$-kernel of $\cF(g_{\bullet})$ in $\Seq(\cA) \colon$
$$\xymatrix{\cF(A_{\bullet}) \ar[rr]^-{\cF(f_{\bullet})} \ar[rrd]_{\sigma} & & \cF(B_{\bullet}) \ar[rr]^-{\cF(g_{\bullet})} & & \cF(C_{\bullet)} \\
& & \cN(\cF(g_{\bullet})) \ar[u]_{n_{\cF(g_{\bullet})}} }$$
Proposition 4.3 in \cite{EVSnail} provides the non trivial part of the proof that the comparison $\sigma$ is a quasi-isomorphism, that is, all
the induced arrows $C(\sigma)_n$ and $K(\sigma)_n$ are isomorphisms. Keeping in mind \ref{TextHomObj}, this means that
$H_n(A_{\bullet}) \simeq H_n(\cF(A_{\bullet})) \simeq \Cok(h_{n+1}^P)$
and we can conclude that the long exact sequence of Corollary \ref{CorLongExSeq} is the usual long homology sequence.
}\end{Text} 

\begin{Text}\label{TextCompChSeq}{\rm
To conclude, let us make more precise the comparison between chain complexes and sequentiable families. Assume that $\cA$ is preadditive
(= enriched in abelian groups). Recall that, in this case, the category $\Ch(\cA)$ is equipped with a structure of nullhomotopies 
$\Theta_{\Ch}$ defined as follows. For a morphism $g_{\bullet} \colon B_{\bullet} \to C_{\bullet},$ a nullhomotopy 
$\varphi_{\bullet} \in \Theta_{\Ch}(g_{\bullet})$ is a family $\{\varphi_n \colon B_n \to C_{n+1}\}_{n \in \mathbb Z}$ of arrows in $\cA$
such that $\varphi_n \cdot d^C_{n+1} + d^B_n \cdot \varphi_{n-1} = g_n$ for all $n \in \mathbb Z$
$$\xymatrix{B_{n+1} \ar[rr]^-{g_{n+1}} \ar[d]_{d^B_{n+1}} & & C_{n+1} \ar[d]^{d^C_{n+1}} \\
B_n \ar[rr]^-{g_n} \ar[d]_{d^B_n} \ar[rru]^{\varphi_n} & & C_n \ar[d]^{d^C_n} \\
B_{n-1} \ar[rr]_-{g_{n-1}} \ar[rru]^{\varphi_{n-1}} & & C_{n-1} }$$
For morphisms $f_{\bullet} \colon A_{\bullet} \to B_{\bullet}$ and $h_{\bullet} \colon C_{\bullet} \to D_{\bullet},$ the nullhomotopy
$f_{\bullet} \circ \varphi_{\bullet} \circ h_{\bullet}$ is defined by $\{f_{n} \cdot \varphi_n \cdot h_{n+1}\}_{n \in \mathbb Z}.$
The structure $\Theta_{\Ch}$ in $\Ch(\cA)$ is not so good as the structure $\Theta_{\Delta}$ in $\Seq(\cA)$ is. Here is why.
\begin{enumerate}
\item The main problem with $\Theta_{\Ch}$ is that the reduced interchange condition \ref{CondRedInter} is not satisfied. To see this,
start with any morphism of complexes $g_{\bullet} \colon B_{\bullet} \to C_{\bullet}$ and any nullhomotopy
$\varphi_{\bullet} \in \Theta_{\Ch}(g_{\bullet}).$ Then, construct the following diagram, where $C_{\bullet}[-1]$ is the $(-1)$-translate 
of $C_{\bullet},$
$$\xymatrix{C_{\bullet} \ar[r] & 0 \ar[r] & C_{\bullet}[-1] }$$
equipped with the nullhomotopy given by the family of identity arrows. Explicitly:
$$\xymatrix{B_{n+1} \ar[rr]^-{g_{n+1}} \ar[d]_{d_{n+1}^B} & & C_{n+1} \ar[rr] \ar[d]^{d_{n+1}^C} & & 0 \ar[rr] \ar[d] & & C_n \ar[d]^{-d_n^C} \\
B_n \ar[rr]_-{g_n} \ar[rru]^{\varphi_n} & & C_n \ar[rr] \ar[rrrru]^<<<<<<<<<<<<{\id} & & 0 \ar[rr] & & C_{n-1} }$$
The reduced interchange condition between the two depicted nullhomotopies gives $\varphi_n \cdot 0 = g_n \cdot \id_{C_n},$ so that 
one would have $g_n = 0$ for all $n \in \mathbb Z.$
\item
Another problem is that, for any complex $C_{\bullet},$ the canonical arrow $\eta_{C_{\bullet}} \colon C_{\bullet} \to \pi_0(C_{\bullet})$
is the identity arrow. This comes from the fact that in the diagram
$$\xymatrix{C_{\bullet} \ar[r] & 0 \ar[r] & C_{\bullet}[-1] }$$
equipped with the already mentioned nullhomotopy, the left part is the $\Theta_{\Ch}$-kernel of the right part and the right part is the 
$\Theta_{\Ch}$-cokernel of the left part. As a consequence, in $(\Ch(\cA),\Theta_{\Ch})$ each object is $S$-global and $S$-proper (in the 
sense of Definition \ref{DefSExact}), and this whatever the class $S$ is.
\end{enumerate} 
}\end{Text}


\begin{Text}\label{TextMorphNull}{\rm
Recall from \cite{VIT23} that a morphism $\cF \colon (\cA,\Theta_{\cA}) \to (\cB,\Theta_{\cB})$ of categories with nullhomotopies  
is a functor $\cF \colon \cA \to \cB$ equipped, for every arrow $g \colon B \to C$ in  $\cA,$ with a map 
$$\cF_g \colon \Theta_{\cA}(g) \longrightarrow \Theta_{\cal B}(\cF(g))$$
such that $\cF_{f \cdot g \cdot h}(f \circ \varphi \circ h) = \cF(f) \circ \cF_g(\varphi) \circ \cF(h)$ for all $f \colon A \to B$ and $h \colon C \to D.$ 
}\end{Text}

\begin{Proposition}\label{PropMorphNull}
If $\cA$ is preadditive, the functor $\cF$ of Proposition \ref{PropSeqAssCh} extends to a morphism of categories with nullhomotopies
$$\cF \colon (\Ch(\cA),\Theta_{\Ch}) \longrightarrow (\Seq(\cA),\Theta_{\Delta})$$
\end{Proposition} 

\begin{proof}
To start, consider a morphism of complexes $g_{\bullet} \colon B_{\bullet} \to C_{\bullet}$ and a nullhomotopy 
$\varphi_{\bullet} \in \Theta_{\Ch}(g_{\bullet}).$ The situation is depicted by the following diagram
$$\xymatrix{ & B_n \ar[rr]^-{g_n} \ar[ddd]_{d^B_n} \ar[ld]_{q^B_n} & & C_n \ar[rd]^{q^C_n} \ar[ddd]^{d^C_n} \\
\Cok(d^B_{n+1}) \ar[rrrr]^-{\overline{\cF g}_n} \ar[d]_{h^{\cF B}_n} & & & & \Cok(d^C_{n+1}) \ar[d]^{h^{\cF C}_n} \\
\Ker(d^B_{n-1}) \ar[rrrr]_-{\underline{\cF g}_n} \ar[rd]_{k^B_{n-1}} & & & & \Ker(d^C_{n-1}) \ar[ld]^{k^C_{n-1}} \\
& B_{n-1} \ar[rr]_{g_{n-1}} \ar[rruuu]^{\varphi_{n-1}} & & C_{n-1} }$$
The nullhomotopy $\cF(\varphi_{\bullet}) \in \Theta_{\Delta}(\cF(g_{\bullet}))$ is defined, in degree $n,$ by the formula
$$\cF(\varphi)_n = k^B_{n-1} \cdot \varphi_{n-1} \cdot q^C_n$$
To check the condition $h^{\cF B}_n \cdot \cF(\varphi)_n = \overline{\cF g}_n,$ compose on the left with the 
epimorphism $q^B_n.$ To check the condition $\cF(\varphi)_n \cdot h^{\cF C}_n = \underline{\cF g}_n,$ compose
on the right with the monomorphism $k^C_{n-1}.$ The rest of the proof is straightforward. 
\end{proof}

\begin{Remark}\label{Rem2Fun}{\rm
Proposition \ref{PropMorphNull} can be improved a bit.
\begin{enumerate}
\item When $\cA$ is preadditive, the category $\Ch(\cA)$ has a 2-categorical structure: an homotopy 
$\varphi_{\bullet} \colon f_{\bullet} \Rightarrow g_{\bullet} \colon B_{\bullet} \rightrightarrows C_{\bullet}$ is an element in 
$\Theta_{\Ch(\cA)}(g_{\bullet}-f_{\bullet}).$ Explicitly, $\varphi_{\bullet}$ is a family $\{\varphi_n \colon B_n \to C_{n+1}\}_{n \in \mathbb Z}$ 
such that $g_n = f_n + \varphi_n \cdot d_{n+1}^C + d_n^B \cdot \varphi_{n-1}$ for all $n \in \mathbb Z.$ A 2-cell 
$[\varphi_{\bullet}] \colon f_{\bullet} \Rightarrow g_{\bullet}$ is a class of homotopies, where two homotopies 
$\varphi_{\bullet},\psi_{\bullet} \colon f_{\bullet} \Rightarrow g_{\bullet}$
are equivalent if there exists a family $\{\alpha_n \colon B_n \to C_{n+2}\}_{n \in \mathbb Z}$ such that
$\varphi_n = \psi_n + \alpha_n \cdot d_{n+2}^C - d_n^B \cdot \alpha_{n-1}$ for all $n \in \mathbb Z.$
\item When $\cA$ is additive and has finite limits, there is an equivalence of categories $\Arr(\cA) \simeq \Grpd(\cA),$ where the latter is
the category of internal groupoids and internal functors in $\cA.$ Since $\Grpd(\cA)$ is a 2-category (2-cells are the internal natural 
transformations), using this equivalence we get a 2-categorical structure on $\Arr(\cA).$ A 2-cell
$\varphi \colon (f,f_0) \Rightarrow (g,g_0) \colon (B,b,B_0) \rightrightarrows (C,c,C_0)$ is an arrow $\varphi \colon B_0 \to C$ such that 
$g = f + b \cdot \varphi$ and $g_0 = f_0 + \varphi \cdot c.$ 
\item Clearly, the 2-categorical structure of $\Arr(\cA)$ can be extended level-wise to a 2-categorical structure for $\Seq(\cA).$
It is easy to see that, in this case, the functor $\cF \colon \Ch(\cA) \to \Seq(\cA)$ of Proposition \ref{PropSeqAssCh} is a 2-functor. 
The argument is essentially the same as in the proof of Proposition \ref{PropMorphNull}, one has just to check that the definition
of $\cF$ on homotopies is compatible with the equivalence relation used in $\Ch(\cA)$ to define 2-cells.

\item Let us end this remark with a point of attention about the terminology introduced in point  \ref{TextCompClassSeqBis}. 
On one hand, the name of quasi-isomorphism is justified by the fact that, by point \ref{TextHomObj}, a morphism in $\Ch(\cA)$ is a 
quasi-isomorphism if and only if its image under the functor $\cF$ is a quasi-isomorphism in $\IsoSeq(\cA).$ On the other hand, by 
analogy with the simpler situation of $\Arr(\cA),$ quasi-isomoprhisms in $\Seq(\cA)$ could be called weak equivalences. Indeed,
under the biequivalence $\Arr(\cA) \simeq \Grpd(\cA)$ recalled above, morphisms $(g,g_0) \colon (B,b,B_0) \to (C;c;C_0)$ such
that $K(g)$ and $C(g_0)$ are isomorphisms correspond to weak equivalences, that is, internal functors which are fully faithful
and essentually surjective.
\end{enumerate}
}\end{Remark}

\hfill

\noindent{\bf Acknowledgements:} We thank Zurab Janelidze for his help with Condition \ref{CondB}
(see Remark \ref{RemZurab}).


\begin{thebibliography}{}
\bibliographystyle{alpha}


\bibitem{BO2} {\sc F. Borceux,}  Handbook of categorical algebra, vol. 2.
Cambridge University Press (1994) xvii+443 pp.

\bibitem{BB} {\sc F. Borceux, D. Bourn,} Mal'cev, protomodular, homological and semi-abelian categories. 
Kluwer Academic Publishers (2004) xiv+479 pp.

\bibitem{DB3x3} {\sc D. Bourn,} $3 \times 3$ lemma and protomodularity, 
{\em Journal of Algebra} 236 (2001) 778--795.

\bibitem{BG} {\sc D. Bourn, M. Gran,} Regular, protomodular and abelian categories.
In: Categorical Foundations, M.C. Pedicchio and W. Tholen Editors, {\em Cambridge University Press} (2004) 169--211.

\bibitem{RB} {\sc R. Brown,} Fibrations of groupoids, {\em Journal of Algebra} 15 (1970) 103--132.

\bibitem{DKV} {\sc J. Duskin, R. Kieboom, E.M. Vitale,} Morphisms of 2-groupoids and
low-dimensional cohomology of crossed modules, {\em Fields Institute Communication Series} 43 (2004) 227--241.

\bibitem{TETVDL} {\sc T. Everaert, T. Van der Linden,} Baer invariants in semi-abelian categories II: homology,
{\em Theory and Applications of Categories} 12 (2004) 195-224.

\bibitem{GZ} {\sc P. Gabriel, M. Zisman,} Calculus of fractions and homotopy theory, 
{\em Springer-Verlag} (1967) x+167 pp.

\bibitem{GR01} {\sc M. Grandis,} A note on exactness and stability in homotopical algebra, 
{\em Theory and Applications of Categories} 9 (2001) 17--42.


\bibitem{PAJ24} {\sc P.-A. Jacqmin,} Surjection-like classes of morphisms,
 {\em Theory and Applications of Categories} 39 (2023) 949--1013

%
\bibitem{JMMVFract} {\sc P.-A. Jacqmin, S. Mantovani, G. Metere, E.M. Vitale,} Bipullbacks of fractions and the snail lemma, 
{\em Journal of Pure and Applied Algebra} 223 (2019) 5147--5162.

\bibitem{ZJ} {\sc Z. Janelidze,} The pointed subobject functor, 3$\times$3 lemmas, and subtractivity of spans,
{\em Theory and Applications of Categories} 23 (2010) 221--242.

\bibitem{ZJEVSnail} {\sc Z. Janelidze, E.M. Vitale,} The snail lemma in a pointed regular category, 
{\em Journal of Pure and Applied Algebra} 221 (2017) 135--143.
%
%
\bibitem{MMV23} {\sc S. Mantovani, M. Messora, E.M. Vitale,} Homotopy torsion theories,
{\em Journal of Pure and Applied Algebra} 228 (2024).
%
\bibitem{MMVSnail} {\sc S. Mantovani, G. Metere, E.M. Vitale,} The snail lemma for internal groupoids, 
{\em Journal of Algebra} 535 (2019) 1--34.

%
%
%
\bibitem{EVSnail} {\sc E.M. Vitale,} The snail lemma, 
{\em Theory and Applications of Categories} 31 (2016) 484--501.
%
\bibitem{VIT23} {\sc E.M. Vitale,} Completion under strong homotopy cokernels, 
{\em Theory and Applications of Categories} 41 (2024) 168--193.

\bibitem{VIT24} {\sc E.M. Vitale,} From abelian categories to 2-abelian bicategories, 
{\em Theory and Applications of Categories} 41 (2024) 1812-1872.

\bibitem{WBook} {\sc C. Weibel,} An introduction to homological algebra,
Cambridge University Press (1994) xiv+450 pp.

\end{thebibliography}
\end{document}